\documentclass[a4paper,12pt,twoside]{article}
\usepackage{amsthm,amssymb,amsmath,amscd}
\usepackage{pb-diagram}
\newtheorem{thm}{Theorem}

\newtheorem{lem}{Lemma}
\newtheorem{prop}{Proposition}
\newtheorem{rem}{Remark}
\newtheorem{conj}{Conjecture}
\theoremstyle{definition}
\newtheorem{defn}{Definition}

\newcommand{\Zset}{\mathbb{Z}}

\newcommand{\bbc}{\mathbb{C}}
\newcommand{\Cset}{\mathbb{C}}

\newcommand{\Nset}{\mathbb{N}}

\newcommand{\ba}{\mathbf{a}}
\newcommand{\bc}{\mathbf{c}}

\def\evarepsilon{{e}}

\def\fp{{\mathfrak p}}
\def\fq{{\mathfrak q}}

\def\cB{{\mathcal B}}

\def\cD{{\mathcal D}}

\def\vol{{\rm vol}}

\def\val{{\rm val}}

\def\Cent{{\rm C}}

\def\cB{{\mathcal B}}
\def\Frob{{\rm Frob}}
\def\cI{{\mathcal I}}

\def\integers{{\mathfrak o}}

\def\SL{{\rm SL}}
\def\GL{{\rm GL}}
\def\SP{{\rm Sp}}

\def\Sp{{\rm Sp}}
\def\Gal{{\rm Gal}}

\def\SO{{\rm SO}}
\def\PGL{{\rm PGL}}

\def\Aut{{\rm Aut}}
\def\alg{{\rm alg}}
\def\diagt{{\rm d}}
\def\Hom{{\rm Hom}}
\def\End{{\rm End}}

\def\Ind{{\rm Ind}}
\def\diag{{\rm diag}}

\def\Mat{{\rm M}}
\def\Nor{{\rm N}}
\def\Irr{{\rm Irr}}

\def\Prim{{\rm Prim}}

\def\cA{{\mathcal A}}
\def\cR{{\mathcal R}}
\def\cH{{\mathcal H}}
\def\cO{{\mathcal O}}
\def\cL{{\mathcal L}}

\def\cW{{\mathcal W}}

\def\fs{{\mathfrak s}}
\def\ft{{\mathfrak t}}

\def\fB{{\mathfrak B}}

\def\integers{{\mathfrak o}}
\def\fR{{\mathfrak R}}
\def\fZ{{\mathfrak Z}}

\def\CC{{\mathbb C}}
\def\fB{{\mathfrak B}}
\def\fR{{\mathfrak R}}

\def\id{{e}}
\def\St{{\rm St}}

\begin{document}
\title{The Hecke algebra of a reductive $p$-adic group: a
geometric conjecture}
\author{Anne-Marie Aubert, Paul Baum and Roger Plymen}
\date{}
\maketitle
 \abstract{Let $\cH(G)$ be the Hecke algebra of a reductive $p$-adic
 group $G$. We formulate a conjecture for the ideals in the Bernstein
 decomposition of $\cH(G)$. The conjecture says that each ideal is
 \emph{geometrically equivalent} to an algebraic variety.   Our
 conjecture is closely related to Lusztig's conjecture on the
 asymptotic Hecke algebra.  We prove our conjecture for $\SL(2)$ and $\GL(n)$.
 We also prove part (1) of the conjecture for the Iwahori ideals of the groups $\PGL(n)$ and $\SO(5)$.}

 \section{Introduction}

 The reciprocity laws in number theory have a long development,
 starting from conjectures of Euler, and including contributions
 of Legendre, Gauss, Dirichlet, Jacobi, Eisenstein, Takagi and Artin. For the details of this development, see
 \cite{Lem}.  The local reciprocity law for a local field $F$, which concerns the
 finite Galois extensions $E/F$ such that $\Gal(E/F)$ is \emph{commutative},
  is stated and proved in \cite [p. 320]{N}.  This local reciprocity law was dramatically
 generalized by Langlands, see \cite{Bbook}. The local Langlands correspondence for $\GL(n)$ is a
 noncommutative generalization of the reciprocity law of local class field theory \cite
 {JR}.   The local Langlands conjectures, and the global Langlands
 conjectures, all involve, inter alia, the representations of
 reductive $p$-adic groups \cite{Bbook}.

 To each reductive $p$-adic group $G$ there is
 associated the Hecke algebra $\cH(G)$, which we now define. Let $K$ be a compact open subgroup of $G$, and define
 $\cH(G//K)$ as the convolution algebra of all complex-valued,
 compactly-supported functions on $G$ such that $f(k_1xk_2) = f(x)$ for
all $k_1$, $k_2$ in $K$.
 The Hecke algebra $\cH(G)$ is then defined as
 \[
 \cH(G): = \bigcup_K \cH(G//K).\]

 The smooth representations of $G$ on a complex vector space $V$ correspond bijectively to the
 nondegenerate representations of $\cH(G)$ on $V$, see \cite[p.2]{B}.

 In this article, we
 consider $\cH(G)$ from the point of view of noncommutative
 (algebraic)
 geometry.

 We recall that the coordinate rings of affine algebraic varieties are
precisely the commutative, unital, finitely generated, reduced $\CC$-algebras,
see \cite[II.1.1]{EH}.

 The Hecke algebra $\cH(G)$ is a non-commutative, non-unital,
non-finitely-generated, non-reduced $\CC$-algebra, and so cannot be the
coordinate ring of an affine algebraic variety.

 The Hecke algebra $\cH(G)$ is non-unital, but it admits \emph{local units},
see \cite[p.2]{B}.

 The algebra $\cH(G)$ admits a canonical decomposition into
ideals, the Bernstein decomposition \cite{B}:
 \[
 \cH(G) = \bigoplus_{\fs \in \fB(G)}
\cH^{\fs}(G).\]

Each ideal $\cH^{\fs}(G)$ is a non-commutative, non-unital,
non-finitely-generated, non-reduced $\CC$-algebra, and so cannot be the
coordinate ring of an affine algebraic variety.

  In section 2, we define the \emph{extended centre} $\widetilde{\fZ}(G)$ of $G$.
  At a crucial point in the construction of the centre $\fZ(G)$ of the
category of smooth representations of $G$, certain
quotients are made: we replace each
 ordinary quotient by the \emph{extended quotient} to create the
 \emph{extended centre}.

 In section 3 we define \emph{morita contexts}, following
 \cite{CDN}.

 In section 4 we prove that each ideal $\cH^{\fs}(G)$ is Morita
 equivalent to a unital $k$-algebra of finite type, where $k$ is the
coordinate ring of a complex affine algebraic variety. We think of
 the ideal $\cH^{\fs}(G)$ as a noncommutative
algebraic variety, and $\cH(G)$ as a noncommutative scheme.

In section 5 we formulate our conjecture. We conjecture that each
ideal $\cH^{\fs}(G)$ is \emph{geometrically equivalent} (in a
sense which we make precise) to the coordinate ring of a complex
affine algebraic variety $X^{\fs}$: \[ \cH^{\fs}(G) \asymp
\cO(X^{\fs}) = \widetilde{\fZ^{\fs}}(G).\] The ring
$\widetilde{\fZ{^\fs}}(G)$ is the $\fs$-factor in the
\emph{extended centre} of $G$.  The ideals $\cH^{\fs}(G)$
therefore qualify as noncommutative algebraic varieties.

We have stripped away the homology and cohomology which play such
a dominant role in \cite{BN}, \cite{BP}, leaving behind three
crucial moves: \emph{Morita equivalence, morphisms which are
spectrum-preserving with respect to filtrations, and deformation
of central character}.
These three moves generate the notion of geometric equivalence.

In section 6 we prove the conjecture for all \emph{generic} points
$\fs \in \fB(G)$.

In section 7 we prove our conjecture for $\SL(2)$.

In section 8 we discuss some general features used in proving the
conjecture for certain examples.

In section 9 we review the asymptotic Hecke algebra of Lusztig.

The asymptotic Hecke algebra $J$ plays a vital role in our
conjecture, as we now proceed to explain. One of the Bernstein
ideals in $\cH(G)$ corresponds to the point $\mathfrak{i} \in
\fB(G)$, where
  $\mathfrak{i}$ is the quotient variety
$\Psi(T)/W_f$. Here, $T$ is a maximal torus in $G$, $\Psi(T)$ is
the complex torus of unramified quasicharacters of $T$, and $W_f$
is the finite Weyl group of $G$.  Let $I$ denote an Iwahori
subgroup of $G$, and define $e$ as follows:
\[
e(x) = \begin{cases}\vol(I)^{-1}&\text{ if $x \in I$,}\cr
0& \text{otherwise.}\end{cases}\]
Then the \emph{Iwahori ideal} is the two-sided ideal generated by
$e$: \[ \cH^{\mathfrak{i}}(G): = \cH(G)e\cH(G).\]  There is a
Morita equivalence $\cH(G) e \cH(G) \sim e \cH(G) e$ and in fact
we have
\[
\cH^{\mathfrak{i}}(G): = \cH(G) e \cH(G) \asymp e \cH(G) e \cong
\cH(G//I) \cong \cH(W, q_F) \asymp J\] where $\cH(W,q_F)$ is an
(extended) affine Hecke algebra based on the (extended) Coxeter
group $W$,  $J$ is the asymptotic Hecke algebra, and $\asymp$
denotes geometric equivalence. Now $J$ admits a decomposition into
finitely many two-sided ideals
\[
J = \bigoplus_\bc J_{\bc}\] labelled by the two-sided cells $\bc$ in $W$.
We therefore have
\[
\cH^{\mathfrak{i}}(G) \asymp \oplus J_{\bc}.\] This canonical
decomposition of $J$ is well-adapted to our conjecture.

In section 10 we prove that
\[J_{\bc_0} \asymp \fZ^{\mathfrak{i}}(G)\] where $\bc_0$ is the lowest
two-sided cell, for any connected $F$-split adjoint simple
$p$-adic group $G$.  We note that $\fZ^{\mathfrak{i}}(G)$ is the
ring of regular functions on the \emph{ordinary quotient}
$\Psi(T)/W_f$.

In section 11 we prove the conjecture for $\GL(n)$.  We establish
that for each point $\fs \in \fB(G)$, we have \[ \cH^{\fs}(\GL(n))
\asymp \widetilde{\fZ^{\fs}}(\GL(n)).\]

In section 12 we prove part (1) of the conjecture for the Iwahori
ideal in $\cH(\PGL(n))$.

In section 13 we prove part (1) of the conjecture for the Iwahori
ideal in $\cH(\SO(5))$. Our proofs depend crucially on Xi's
affirmation, in certain special cases, of Lusztig's conjecture on
the asymptotic Hecke algebra $J$ (see \cite[\S 10]{LCellsIV}).

In section 14 we discuss some consequences of the conjecture.

We thank the referee for his many detailed and constructive
comments, which forced us thoroughly to revise the article. We
also thank Nigel Higson, Ralf Meyer and Victor Nistor for valuable
discussions, and Gene Abrams for an exchange of emails concerning
rings with local units.

\section{The extended centre}\label{centre}

Let $G$ be the set of rational points of a reductive group defined
over a local nonarchimedean field $F$, and let $\fR(G)$ denote the
category of smooth $G$-modules. Let $(L, \sigma)$ denote a
\emph{cuspidal pair}: $L$ is a Levi subgroup of $G$ and $\sigma$
is an irreducible supercuspidal representation of $L$. The group
$\Psi(L)$ of unramified quasicharacters of $L$ has the structure
of a complex torus.

We write $[L,\sigma]_G$ for the equivalence class of $(L,\sigma)$ and
$\fB(G)$ for the set of equivalence classes, where the equivalence relation
is defined by $(L,\sigma)\sim(L',\sigma')$ if $gLg^{-1}=L'$ and
${}^g\sigma\simeq\nu'\sigma'$, for some $g\in G$ and some $\nu'\in\Psi(L')$.
For $\fs=[L,\sigma]_G$, let $\fR^\fs(G)$ denote the full subcategory of
$\fR(G)$ whose objects are the representations $\Pi$ such that each
irreducible subquotient of $\Pi$ is a subquotient of a parabolically induced
representation $\iota_P^G(\nu\sigma)$ where $P$ is a parabolic subgroup of $G$
with Levi subgroup $L$ and $\nu\in\Psi(L)$.
The action (by conjugation) of $\Nor_G(L)$ on $L$ induces an action of
$W(L)=\Nor_G(L)/L$ on $\fB(L)$. Let $W_\ft$ denote the stabilizer of
$\ft=[L,\sigma]_L$ in $W(L)$. Thus $W_\ft=N_\ft/L$ where
\[N_\ft=\left\{n\in\Nor_G(L)\,:\,{}^n\sigma\simeq\nu\sigma, \;\text{ for
some $\nu\in \Psi(L)$}\right\}.\]
It acts (via conjugation) on $\Irr^\ft\, L$, the set of isomorphism classes
of irreducible objects in $\fR^\ft(L)$.

Let $\Omega(G)$ denote the set of $G$-conjugacy classes of
cuspidal pairs $(L,\sigma)$. The groups $\Psi(L)$ create orbits in
$\Omega(G)$. Each orbit is of the form $D_{\sigma}/W_{\ft}$ where
$D_{\sigma}=\Irr^\ft L$ is a complex torus.

We have
\[
\Omega(G) = \bigsqcup \, D_{\sigma}/W_{\ft}.
\]

Let $\fZ(G)$ denote the \emph{centre} of the category $\fR(G)$.
The centre of an abelian category (with a small skeleton) is the
endomorphism ring of the identity functor. An element $z$ of the
centre assigns to each object $A$ in $\fR(G)$ a morphism $z(A)$
such that
\[f\cdot z(A) = z(B)\cdot f\] for each morphism $f \in \Hom(A,B)$.

According to Bernstein's theorem \cite{B} we have the explicit
decomposition of $\fR(G)$:
\[\fR(G) = \prod_{\fs\in\fB(G)} \fR^{\fs}(G).\]
We also have
\[
\fZ(G) \cong \prod \fZ^{\fs}\] where \[ \fZ^{\fs}(G) =
\cO(D_{\sigma}/W_{\ft})\]
is the centre of the category $\fR^{\fs}(G)$.

Let the finite group $\Gamma$ act on the space $X$. We define, as
in \cite{BC},
\[
\widetilde{X}: = \{(\gamma,x): \gamma x = x\} \subset \Gamma
\times X\] and define the $\Gamma$-action on $\widetilde{X}$ as
follows:
\[
\gamma_1(\gamma,x): = (\gamma_1 \gamma \gamma_1^{-1}, \gamma_1
x).\] The \emph{extended quotient} of $X$ by $\Gamma$ is defined
to be the ordinary quotient $\widetilde{X}/\Gamma$. If $\Gamma$ acts
freely, then we have
\[
\widetilde{X} = \{(1,x): x \in X\} \cong X
\]
and, in this case, $\widetilde{X}/\Gamma = X/\Gamma$.

We will write
\[
\widetilde{\fZ}^{\fs}(G): =
\mathcal{O}(\widetilde{D_{\sigma}}/W_{\ft}).\]

We now form the \emph{extended centre}
\[\widetilde{\fZ}(G): = \prod_{\fs \in \fB(G)} \mathcal{O}(\widetilde{D_{\sigma}}/W_{\ft}).\]

We will write
\[
k: = \cO(D_{\sigma}/W_{\ft}) = \fZ^{\fs}(G).\]

\section{Morita contexts and the central character}

A fundamental theorem of Morita says that the categories of
modules over two rings with identity $R$ and $S$ are equivalent if
and only if there exists a strict Morita context connecting $R$
and $S$.

A Morita context connecting $R$ and $S$ is a datum \[(R,S,_RM_S,
_SN_R,\phi,\psi)\] where $M$ is an $R$-$S$-bimodule, $N$ is an
$S$-$R$-bimodule, $\phi: M\otimes_SN \to R$ is a morphism of
$R$-$R$-bimodules and $\psi:N \otimes_RM \to S$ is a morphism of
$S$-$S$-bimodules such that
\[
\phi(m \otimes n)m' = m\psi(n \otimes m')\]
 \[\psi(n \otimes m)n'
= n\psi(m \otimes n')\] for any $m,m' \in M, n,n'\in N.$  A Morita
context is \emph{strict} if both maps $\phi,\psi$ are
isomorphisms, see \cite{CDN}.

Let $R$ be a ring with local units, i.e. for any finite subset $X$
of $R$, there exists an idempotent element $e \in R$ such that $ex
= xe = x$ for any $x \in X$.  An $R$-module $M$ is unital if $RM =
M.$  Unital modules are called non-degenerate in \cite{B}.

Let $R$ and $S$ be two rings with local units, and $R-MOD$ and
$S-MOD$ be the associated categories of unital modules.  A Morita
context for $R$ and $S$ is a datum $(R,S,_RM_S, _SN_R,\phi,\psi)$
with the condition that $M$ and $N$ are unital modules to the left
and to the right.  If the Morita context is strict, then we obtain
by \cite[Theorem 4.3]{CDN} an equivalence of categories $R-MOD$
and $S-MOD$.

If $R$ and $S$ are rings with local units which are connected by a
strict Morita context, then we will say that $R$ and $S$ are
\emph{Morita equivalent} and write
\[
R \sim_{morita} S.\]

We will use repeatedly the following elementary lemmas.

\begin{lem}\label{matrix}  Let $R$ be a ring with local units.
Let $M_n(R)$ denote $n \times n$
matrices over $R$.  Then we have
\[
M_n(R) \sim_{morita} R.\]
\end{lem}

\begin{proof}  Let $M_{i \times j}(R)$ denote $i \times j$
matrices over $R$.  Then we have a strict Morita context
\[
(R, M_{n\times n}(R), M_{1 \times n}(R), M_{n \times 1}(R), \phi,
\psi)\] where $\phi, \psi$ denote matrix multiplication.
\end{proof}

\begin{lem}\label{local}  Let $R$ be a ring with local units. Let $e$ be an idempotent in $R$.
Then we have
\[
ReR \sim_{morita} eRe.\]
\end{lem}

\begin{proof}  Given $r \in R$, there is an idempotent $e \in R$ such that $r = er \in
R^2$ so that $R \subset R^2 \subset R$. A ring $R$ with local
units is \emph{idempotent}: $R = R^2$.  This creates a strict
Morita context \[(eRe, ReR, eR, Re, \phi, \psi)\]  where $\phi,
\psi$ are the obvious multiplication maps in $R$.  We now check
identities such as $(eR)(Re) = eRe, (Re)(eR) = ReR$.
\end{proof}

Let $X$ be a complex affine algebraic variety, let
$\mathcal{O}(X)$ denote the algebra of regular functions on $X$.
Let\[k: = \mathcal{O}(X).\]

A $k$-algebra $A$ is a $\Cset$-algebra which is also a $k$-module
such that \[1_k\cdot a = a, \quad \omega (a_1a_2) = (\omega
a_1)a_2 = a_1(\omega a_2), \quad \omega(\lambda a) =
\lambda(\omega a) = (\lambda \omega)a\] for all $a,a_1, a_2 \in A,
\omega \in k, \lambda \in \Cset$.

If $A$ has a unit then a $k$-algebra is a $\Cset$-algebra with a
given unital homomorphism of $\Cset$-algebras from $k$ to the
centre of $A$. If $A$ does not have a unit, then a $k$-algebra is
a $\Cset$-algebra with a given unital homomorphism of
$\Cset$-algebras from $k$ to the centre of the multiplier algebra
$\mathcal{M}(A)$ of $A$, see \cite{B}.

A $k$-algebra $A$ is of \emph{finite type} if, as a $k$-module,
$A$ is finitely generated.

If $A,B$ are $k$-algebras, then a morphism of $k$-algebras is a
morphism of $\Cset$-algebras which is also a morphism of
$k$-modules.

If $A$ is a $k$-algebra and $Y$ is a unital $A$-module then the
action of $A$ on $Y$ extends uniquely to an action of
$\mathcal{M}(A)$ on $Y$ such that $Y$ is a unital
$\mathcal{M}(A)$-module \cite{B}. In this way the category $A-MOD$
of unital $A$-modules is equivalent to the category
$\mathcal{M}(A)-MOD$ of unital $\mathcal{M}(A)$-modules.  Since
the $k$-algebra structure for $A$ can be viewed as a unital
homomorphism of $\Cset$-algebras from $k$ to the centre of
$\mathcal{M}(A)$, it now follows that given any unital $A$-module
$Y$, $Y$ is canonically a $k$-module with the following
compatibility between the $A$-action and the $k$-action:
\[1_k\cdot y = y, \quad \omega (ay) = (\omega
a)y = a(\omega y), \quad \omega(\lambda y) = \lambda(\omega y) =
(\lambda \omega)y
\] for all $a \in A, y \in Y, \omega \in k, \lambda \in \Cset$.

If $A,B$ are $k$-algebras then a strict Morita context (in the
sense of $k$-algebras) connecting $A$ and $B$ is a $6$-tuple
\[
(A,B,M,N,\phi,\psi)\] which is a strict Morita context connecting
the rings $A,B$. We require, in addition, that
\[\omega y = y \omega \] for all $y \in M, \omega
\in k$.  Similarly for $N$.  When this is satisfied, we will say
that the $k$-algebras $A$ and $B$ are \emph{Morita equivalent} and
write
\[
A \sim_{morita}B.\]

By \emph{central character} or \emph{infinitesimal character} we
mean the following.   Let $M$ be a simple $A$-module.  Schur's
lemma implies that each $\theta \in k$ acts on $M$ via a complex
number $\lambda(\theta)$.   Then $\theta \mapsto \lambda(\theta)$
is a morphism of $\Cset$-algebras $k \to \Cset$ and is therefore
given by evaluation at a $\Cset$-point of $X$.   The map $\Irr(A)
\to X$ so obtained is the central character (or infinitesimal
character). The notation for the infinitesimal character is as
follows:
\[inf.\,ch.\,: \Irr(A) \longrightarrow X.\]

If $A$ and $B$ are $k$-algebras connected by a strict Morita
context then we have a commutative diagram in which the top
horizontal arrow is bijective:

\[
\begin{CD}
\Irr(A) @>>> \Irr(B)\\
@V inf.ch.VV                @VV inf.ch. V\\
X @>>id > X \end{CD}\]

The algebras that occur in this paper have the property that
\[
\Prim(A) = \Irr(A).\] In particular, each Bernstein ideal
$\cH^{\fs}(G)$ has this property and any $k$-algebra of finite
type has this property, see \cite{BN}.

\section{A Morita equivalence}

Let $\cH = \cH(G)$ denote the Hecke algebra of $G$.  Note that
$\cH(G)$ admits a set $E$ of local units.   For let $K$ be a
compact open subgroup of $G$ and define
\[
e_K(x) = \begin{cases}\vol(K)^{-1}&\text{ if $x \in K$,}\cr 0&
\text{otherwise.}\end{cases}\]

Given a finite set $X \subset \cH(G)$, we choose $K$ sufficiently
small. Then we have $e_Kx = x = xe_K$  for all $x \in X$.  It
follows that $\cH$ is an idempotent algebra: $\cH = \cH^2$.

We have the Bernstein decomposition
\[
\mathcal{H}(G) = \bigoplus_{\fs \in \fB(G)} \cH^{\fs}(G)
\]
of the Hecke algebra $\cH(G)$ into two-sided ideals.


\begin{lem} \label{unitlemma}  Each Bernstein ideal $\cH^{\fs}(G)$  admits a set of local
units.
\end{lem}

\begin{proof} Define
\[E^{\fs}: = E \cap \cH^{\fs}(G).\]
Then $E^{\fs}$ is a set of local units for $\cH^{\fs}(G)$.
\end{proof}

We recall the notation from section~\ref{centre}:
\[\fs \in \fB(G),\quad \fs = [L,\sigma]_G,\quad D_{\sigma} =
\Psi(L)/\mathcal{G}.\]

\begin{thm}\label{one}
Let $\fs \in \fB(G), k = \mathcal{O}(D_{\sigma})$. The ideal
$\cH^{\fs}(G)$ is a $k$-algebra Morita equivalent to a unital
$k$-algebra of finite type.  If $W_{\ft} = \{1\}$ then
$\cH^{\fs}(G)$ is Morita equivalent to $k$.
\end{thm}
\begin{proof}  Let $\fs = [L,\sigma]_G$.  We will write
\[
\cH(G) = \cH = \cH^{\fs} \oplus \cH'\] where $\cH'$ denotes the
sum of all $\cH^{\ft}$ with $\ft \in \fB(G)$ and $\ft \neq \fs$.
It follows from \cite[(3.7)]{B} that there is a compact open
subgroup $K$ of $G$ with the property that $V^K \neq 0$ for every
irreducible representation $(\pi,V)$ with $\mathfrak{I}(\pi) =
\fs$.  We will write $e_K = e + e'$ with $e \in \cH^{\fs}, e' \in
\cH'$.  Both $e$ and $e'$ are idempotent and we have $ee' = 0 =
e'e$.

Given $h \in \cH$, we will write $h = h^{\fs} + h'$ with $h^{\fs}
\in \cH^{\fs}, h' \in \cH'$. By \cite[3.1,\,3.4 -- 3.6]{BKtyp}, we
have $\cH^{\fs}(G) \cong \cH e \cH$.  We note that this follows
from the general considerations in \cite[\S 3]{BKtyp}, and does
\emph{not} use the existence or construction of types.  The ideal
$\cH^{\fs}$ is the idempotent two-sided ideal generated by $e$. By
Lemma~\ref{local} we have
\[
\cH^{\fs} = \cH e \cH \sim_{morita} e\cH e.\]  Since $eh' = he' =
0, e'h^{\fs} = h^{\fs}e' = 0$ we also have
\[
e\cH e = e\cH^{\fs} e = e_K\cH^{\fs}e_K.\] It follows that
\[
\cH^{\fs} \sim_{morita} e_K\cH^{\fs} e_K.\]

Let $B$ be the unital algebra defined as
follows:
\[
B: = e_K\cH^{\fs}(G)e_K.\]

We will use the notation in \cite[p.80]{BNotes}. Let $\Omega =
D_{\sigma}/W_{\ft}$ and let $\Pi: = \Pi(\Omega)^K, \Lambda: =
\End\,\Pi$. There is a strict Morita context
\[
(B, \Lambda, \Pi, \Pi^*, \phi, \psi).\]

We therefore have
\[
\cH^{\fs}(G) \sim_{morita} B \sim_{morita} \Lambda.\]

The intertwining operators $A_w$ generate $\Lambda$ over
$\Lambda(D_{\sigma})$, see \cite[p.81]{BNotes}.  We also have
\cite[p. 73]{BNotes}
\[
\Lambda(D_{\sigma}) = \mathcal{O}(\Psi(L)) \rtimes \mathcal{G}\]
where the finite abelian group $\mathcal{G}$ is defined as
follows:
\[ \mathcal{G}: = \{\psi \in \Psi(G): \psi \sigma \cong
\sigma\}.\] Note that $\mathcal{G}$ acts freely on the complex
torus $\Psi(L)$.

The $k$-algebras $\mathcal{O}(\Psi(L)) \rtimes \mathcal{G}$ and
$\mathcal{O}(\Psi(L)/\mathcal{G})$ are connected by the following
strict Morita context:
\[
(\mathcal{O}(\Psi(L)) \rtimes \mathcal{G},
\mathcal{O}(\Psi(L)/\mathcal{G}), \mathcal{O}(\Psi(L)),
\mathcal{O}(\Psi(L)), \phi, \psi).\]

We conclude that
 \[\Lambda(D_{\sigma}) = \mathcal{O}(\Psi(L))
\rtimes \mathcal{G} \sim_{morita} \mathcal{O}(\Psi(L)/\mathcal{G})
= \mathcal{O}(D_{\sigma}) = k.\]  Therefore, $\Lambda$ is finitely
generated as a $k$-module.  Note that  \[
Centre\,\Lambda(D_{\sigma}) \cong k.\]

If $W_{\ft} = \{1\}$ then there are no intertwining operators, so
that $\Lambda = \Lambda(D_{\sigma})$. In that case we have
\[
\cH^{\fs}(G) \sim_{morita} k.\]

\end{proof}

\section{The conjecture}

Let $k$ be the coordinate ring of a complex affine algebraic
variety $X$, $k = \mathcal{O}(X)$.  Let $A$ be an associative
$\Cset$-algebra which is also a $k$-algebra.  We work with the
collection of all $k$-algebras $A$ which are countably generated.
As a $\Cset$-vector space, $A$ admits a finite or countable basis.


We will define an equivalence relation, called \emph{geometric
equivalence}, on the collection of such algebras $A$.  This
equivalence relation will be denoted $\asymp$.

(1) \textsc{Morita equivalence of $k$-algebras with local units}.
Let $A$ and $B$ be $k$-algebras, each with a countable set of
local units.
If $A$ and $B$ are connected by a strict Morita context, then $A
\asymp B$.    Periodic cyclic homology is preserved, see
\cite[Theorem 1]{C}.

(2)  \textsc{Spectrum preserving morphisms with respect to
filtrations of $k$-algebras of finite type}, as in \cite{BN}. Such
filtrations are automatically finite, since $k$ is Noetherian.

A morphism $\phi\colon A \to B$ of $k$-algebras of finite type
is called
\begin{itemize}
\item \emph{spectrum preserving} if, for each primitive ideal
$\fq$ of $B$, there exists a unique primitive ideal $\fp$ of $A$
containing $\phi^{-1}(\fq)$, and the resulting map $\fq\mapsto\fp$
is a bijection from $\Prim(B)$ onto $\Prim(A)$; \item
\emph{spectrum preserving with respect to filtrations} if there
exist increasing filtrations by ideals
\[
(0)=I_0 \subset I_1 \subset I_2 \subset \cdots \subset I_{r-1}
\subset I_r \subset A\]
\[
(0)=J_0 \subset J_1 \subset J_2 \subset \cdots \subset J_{r-1}
\subset J_r \subset B,\] such that, for all $j$, we have
$\phi(I_j)\subset J_j$ and the induced morphism
\[
\phi_{*}: I_j/I_{j-1} \to J_j/J_{j-1}\] is spectrum preserving. If
$A,B$ are $k$-algebras of finite type such that there exists a
morphism $\phi : A \to B$ of $k$-algebras which is spectrum
preserving with respect to filtrations then $A \asymp B$.
\end{itemize}
(3) \textsc{Deformation of central character}. Let $A$ be a unital
algebra over the complex numbers. Form the algebra $A[t,t^{-1}]$
of Laurent polynomials with coefficients in $A$. If $q$ is a
non-zero complex number, then we have the evaluation-at-$q$ map of
algebras
\[A[t,t^{-1}] \longrightarrow A\]
which sends a Laurent polynomial $P(t)$ to $P(q)$. Suppose now
that $A[t,t^{-1}]$ has been given the structure of a $k$-algebra
i.e. we are given a unital map of algebras over the complex
numbers from $k$ to the centre of $A[t,t^{-1}]$. We assume that
for any non-zero complex number $q$ the composed map
\[
k \longrightarrow A[t,t^{-1}] \longrightarrow A\] where the second
arrow is the above evaluation-at-$q$ map makes $A$ into a
$k$-algebra of finite type. For $q$ a non-zero complex number,
denote the finite type k-algebra so obtained by $A(q)$. Then we
decree that if $q_1$ and $q_2$ are any two non-zero complex
numbers, $A(q_1)$ is equivalent to $A(q_2)$.

We fix $k$. The first two moves preserve the central character.
This third move allows us to algebraically deform the central
character.

Let $\asymp$ be the equivalence relation generated by (1), (2), (3);
we say that $A$ and $B$ are \emph{geometrically equivalent} if $A \asymp B$.

Since each move induces an isomorphism in periodic cyclic homology
\cite{BN} \cite{C}, we have
\[
A \asymp B \Longrightarrow HP_*(A) \cong HP_*(B).\]

\smallskip

In order to formulate our conjecture, we need to review certain
results and definitions.

The primitive ideal space of $\widetilde{\fZ^{\fs}}(G)$ is the set
of $\mathbb{C}$-points of the variety
$\widetilde{D_{\sigma}}/W_{\ft}$ in the Zariski topology.

 We have an isomorphism
\[HP_*(\cO(\widetilde{D_{\sigma}}/W_{\ft})) \cong
H^*(\widetilde{D_{\sigma}}/W_{\ft};\Cset).\]  This is a special
case of the Feigin-Tsygan theorem; for a proof of this theorem
which proceeds by reduction to the case of smooth varieties, see
\cite{KNS}.

Let $E_{\sigma}$ be the maximal compact subgroup of the complex
torus $D_{\sigma}$, so that $E_{\sigma}$ is a compact torus.

Let $\Prim^{t}\,\cH^{\fs}(G)$ denote the set of primitive ideals
attached to tempered, simple $\cH^{\fs}(G)$-modules.

\begin{conj} Let $\fs \in \fB(G)$. Then

{\rm (1)} $\cH^{\fs}(G)$ is geometrically equivalent to the commutative
algebra $\widetilde{\fZ^{\fs}}(G)$:
\[ \cH^{\fs}(G)\asymp \widetilde{\fZ^{\fs}}(G).\]

{\rm (2)} The resulting bijection of primitive ideal spaces
\[
\Prim \, \cH^{\fs}(G) \longleftrightarrow
\widetilde{D_{\sigma}}/W_{\ft}
\]
restricts to give a bijection
\[\Prim^{t}\,\cH^{\fs}(G)
\longleftrightarrow \widetilde{E_{\sigma}}/W_{\ft}.\]

\end{conj}

\begin{rem}{\rm The geometric equivalence of part (1) induces an isomorphism
\[HP_*(\cH^{\fs}(G)) \cong H^*(\widetilde{D_\sigma}/W_{\ft};\Cset).
\]
}
\end{rem}

\begin{rem}{\rm

The referee has posed a very interesting question: if $A,B$ are
geometrically equivalent, what does this imply about the
categories $A-MOD$ and $B-MOD$ of unital modules?

It seems likely that for each ideal $\cH^{\fs}$ the category
$\cH^{\fs}-MOD$ will have some resemblance to
$\widetilde{\fZ}^{\fs}-MOD$.  If the conjecture is true, then
$\Prim\,\cH^{\fs}$ is in bijection with
$\widetilde{D}_{\sigma}/W_{\ft}$.  However, as the referee has
indicated, there may be further resemblances between
$\cH^{\fs}-MOD$ and $\widetilde{\fZ}^{\fs}-MOD$. }
\end{rem}

\section{Generic points in the Bernstein spectrum}

We begin with a definition.

\begin{defn}  The point $\fs \in \fB(G)$
is \emph{generic} if $W_\ft = \{1\}$.
\end{defn}

For example, let $\fs = [G,\sigma]_G$ with $\sigma$ an irreducible
supercuspidal representation of $G$. Then $\fs$ is a generic point
in $\fB(G)$.  For a second example, let $\fs = [\GL(2)\times
\GL(2), \sigma_1 \otimes \sigma_2]_{\GL(4)}$ with $\sigma_1$ not
equivalent to $\sigma_2$ (after unramified twist). Then $\fs$ is a
generic point in $\fB(\GL(4))$.

\begin{thm} \label{regular}  The conjecture is true if $\fs$ is a generic point in
$\fB(G)$.
\end{thm}
\begin{proof}  Part (1).  This is immediate from Theorem~\ref{one}.
We conclude that
\[
\cH^{\fs}(G) \asymp \mathcal{O}(D_{\sigma}) = \fZ^{\fs}(G) =
\widetilde{\fZ}^{\fs}(G).\]



Part (2).   Let $\mathcal{C}(G)$ denote the Harish-Chandra
Schwartz algebra of $G$, see \cite{Wald}.
We choose $e, e_K$ as in the proof of Theorem~\ref{one}.  Let $\fs
\in \fB(G)$ and let $\mathcal{C}^{\fs}(G)$ be the corresponding
Bernstein ideal in $\mathcal{C}(G)$.  As in the proof of
Theorem~\ref{one}, $\mathcal{C}^{\fs}(G)$ is the two-sided ideal
generated by $e$. We have $\mathcal{C}^{\fs}(G) = \mathcal{C}(G) e
\mathcal{C}(G)$. By Lemma~\ref{local}, we have
\[
\mathcal{C}^{\fs}(G) = \mathcal{C}(G)e\mathcal{C}(G) \sim_{morita}
e\mathcal{C}(G) e = e_K\mathcal{C}(G)e_K.\]

According to Mischenko's
theorem \cite{Mis}, the Fourier transform induces an isomorphism
of unital Fr\'{e}chet algebras:
\[
e_K\mathcal{C}(G)e_K \cong C^{\infty}(E_{\sigma}, \End \,E^K).\]

We also have \[C^{\infty}(E_{\sigma}, \End \,E^K) \cong
M_n(C^{\infty}(E_{\sigma}))\]with $n = \dim_{\Cset}(E^K)$.  By
Lemma~\ref{matrix} we have
\[
\mathcal{C}^{\fs}(G) \sim_{morita} C^{\infty}(E_{\sigma}).\]

We now exploit the liminality of the reductive group $G$, and take
the primitive ideal space of each side. We conclude that there is
a bijection
\[
\Prim^t \cH^{\fs}(G) \longleftrightarrow E_{\sigma}.\]

\end{proof}

\section{The Hecke algebra of $\SL(2)$}

Let $G = \SL(2) = \SL(2,F)$, a group not of adjoint type.  Let $W$
be the Coxeter group with $2$ generators:
\[
W = \langle s_1,s_2 \rangle = \mathbb{Z} \rtimes W_f\] where $W_f
= \mathbb{Z}/2\mathbb{Z}$. Then $W$ is the infinite dihedral
group. It has the property that
\[
\cH(W,q_F) = \cH(\SL(2)//I).
\]
There is a unique isomorphism of $\Cset$-algebras between
$\cH(W,q_F)$ and $\Cset[W]$ such that
\[
T_{s_1} \mapsto \frac{q_F + 1}{2} \cdot s_1 + \frac{q_F
-1}{2},\quad T_{s_2} \mapsto \frac{q_F+1}{2} \cdot s_2 +
\frac{q_F-1}{2},\]
where $T_{s_i}$ is the element of $\cH(W,q_F)$ corresponding to $s_i$.
We note also that

\[
\Cset[\Zset \rtimes \Zset/2\Zset] \cong \mathbb{M} \rtimes
\Zset/2\Zset\] where
\[\mathbb{M}: = \Cset[t,t^{-1}]\]
denotes the $\mathbb{Z}/2\mathbb{Z}$-graded algebra of Laurent
polynomials in one indeterminate $t$.   Let $\alpha$ denote the
generator of $\Zset/2\Zset$.  The group $\Zset/2\Zset$ acts as
automorphism of $\mathbb{M}$, with $\alpha(t) = t^{-1}$.  We
define
\[
\mathbb{L}: = \{P\in \mathbb{M}: \alpha(P) = P\}\] as the algebra
of balanced Laurent polynomials.   We will write
\[
\mathbb{L}^*: = \{P \in \mathbb{M}: \alpha(P) = - P\}.\] Then
$\mathbb{L}^*$ is a free of rank $1$ module over $\mathbb{L}$,
with generator $t - t^{-1}$.  We will refer to the elements of
$\mathbb{L}^*$ as \emph{anti-balanced} Laurent polynomials.
\begin{lem} \label{secondlemma} We have
\[\mathbb{M} \rtimes \Zset/2\Zset \asymp \Cset^2 \oplus
\mathbb{L}.\]
\end{lem}
\begin{proof}

We will realize the crossed product as follows:
 \[
 \mathbb{M} \rtimes \Zset/2\Zset = \{f \in
 \mathcal{O}(\Cset^{\times}, M_2(\Cset))\; :\; f(z^{-1}) = \mathfrak{a}
 \cdot f(z) \cdot \mathfrak{a}^{-1}\}\]
 where
 \[\mathfrak{a} =
\left(\begin{matrix}1&0\\
0&-1\end{matrix}\right)\] We then have
\[\mathbb{M} \rtimes \Zset/2\Zset =
\left(\begin{matrix}\mathbb{L}&\mathbb{L}^*\\
\mathbb{L}^*&\mathbb{L}\end{matrix}\right)\]  There is an algebra
map
\[
\left(\begin{matrix}\mathbb{L}&\mathbb{L}^*\\
\mathbb{L}^*&\mathbb{L}\end{matrix}\right) \to
\left(\begin{matrix}\mathbb{L}&\mathbb{L}\\
\mathbb{L}&\mathbb{L}\end{matrix}\right)
\]
as follows: \[x_{11} \mapsto x_{11},\; x_{22} \mapsto x_{22},\;
x_{12} \mapsto x_{12}(t - t^{-1}),\; x_{21} \mapsto x_{21}(t -
t^{-1})^{-1}.\] This map, combined with evaluation of $x_{22}$ at
$t = 1$, $t = -1$, creates an algebra map\[ \mathbb{M} \rtimes
\Zset/2\Zset \to M_2(\mathbb{L}) \oplus \Cset \oplus \Cset.\] This
map is spectrum-preserving with respect to the following
filtrations:
\[
0 \subset M_2(\mathbb{L}) \subset M_2(\mathbb{L}) \oplus \Cset
\oplus \Cset\]
\[
0 \subset I \subset \mathbb{M} \rtimes \mathbb{Z}/2\mathbb{Z}
\]
where $I$ is the ideal of $\mathbb{M} \rtimes
\mathbb{Z}/2\mathbb{Z}$ defined by the conditions \[x_{22}(1) = 0
= x_{22}(-1).\] We therefore have
\[ \mathbb{M} \rtimes \Zset/2\Zset \asymp M_2(\mathbb{L}) \oplus
\Cset^2 \asymp \mathbb{L} \oplus \Cset^2.\]

It is worth noting that $\mathbb{M} \rtimes
\mathbb{Z}/2\mathbb{Z}$ is \emph{not} Morita equivalent to
$\mathbb{L} \oplus \Cset^2$.  For the primitive ideal space
$\Prim(\mathbb{M} \rtimes \mathbb{Z}/2\mathbb{Z})$ is connected,
whereas the primitive ideal space $\Prim(\mathbb{L} \oplus
\Cset^2)$ is disconnected (it has $3$ connected components).
\end{proof}

\begin{thm}~\label{elliptic}
The conjecture is true for $\SL(2,F)$.
\end{thm}

\begin{proof}  For the Iwahori ideal $\cH^{\mathfrak{i}}(\SL(2))$
we have \[\cH^{\mathfrak{i}}(G) \asymp \cH(W,q_F) \asymp
\Cset[W].\]

Let $\chi$ be a unitary character of $F^{\times}$ of exact order
two. Let $G=\SL(2,F)$ and let $\mathfrak{j}= \mathfrak{j}(\lambda)
= [T,\lambda]_G\in\fB(G)$ with $\lambda$ is defined as follows:
$$\lambda\colon\left(\begin{matrix}x&0\\
 0&x^{-1}\end{matrix}\right) \mapsto \chi(x).
$$
Let $c$ be the least integer $n\ge 1$ such that
$1+\fp_F^n\subset\ker(\chi)$. We set
$$J_\chi:=\left(\begin{matrix}\integers_F^\times&\fp_F^{[c/2]}\cr
\fp_F^{[(c+1)/2]}&\integers_F^\times\cr
\end{matrix}\right)\cap\SL(2,F).$$
Let $\tau_\chi$ denote the restriction of $\lambda$ to the compact torus
$$\left\{\left(\begin{matrix}x&0\\
 0&x^{-1}\end{matrix}\right)\,:\,x\in\integers_F^\times\right\}.$$
Then $(J_\chi,\tau_\chi)$ is an $\fs$-type in $\SL(2,F)$ and (for
instance, as a special case of \cite[Theorem~11.1]{GR2}) the Hecke
algebra $\cH(G,\tau_\chi)$ is isomorphic to $\cH(W,q_F)$.  We have
\[
\cH^{\mathfrak{j}}(G) \asymp \cH(G,\tau_{\chi}) \cong \cH(W,q_F)
\asymp \Cset[W].\]

To summarize: if $\fs = \mathfrak{i}$ or $\mathfrak{j}(\lambda)$
then we have
\[
\cH^{\fs}(G) \asymp \mathbb{M} \rtimes \Zset/2\Zset.\]

We have $W_\ft=W_f = \Zset/2\Zset$. Let $\Omega$ be the variety
which corresponds to $\mathfrak{j} \in \fB(G)$, so that $\Omega =
D/W_f$ with $D$ a complex torus of dimension $1$. Each unramified
quasicharacter of the maximal torus $T \subset \SL(2)$ is given by
\[\left(\begin{matrix}x&0\\
 0&x^{-1}\end{matrix}\right) \mapsto s^{\val(x)}\]
 with $z \in \Cset^{\times}$.  Let $\mathbb{V}(zw-1)$ denote the algebraic curve
 in $\Cset^2$ defined by the equation $zw-1 = 0$.  This algebraic curve is a hyperbola.
 The map
 $\Cset^{\times} \to \mathbb{V}(zw-1), \,z \mapsto (z,z^{-1})$ defines the
 structure of algebraic curve on $\Cset^{\times}$.
 The generator of $W_f$ sends a point $z$ on this curve to $z^{-1}$ and there are two
 fixed points, $1$ and $-1$. The coordinate algebra of the quotient curve is given by
 \[
 \Cset[D/W_f] = \mathbb{L}\] and the coordinate algebra of the
 extended quotient is given by
 \[
 \Cset[\widetilde{D}/W_f] = \mathbb{L} \oplus \Cset \oplus
 \Cset.\]
 Then it follows from Lemma~\ref{secondlemma} that
\[
\cH^{\fs}(\SL(2)) \asymp \widetilde{\fZ}^{\fs}(\SL(2))\] if $\fs =
\mathfrak{i}$ or $\mathfrak{j}(\lambda)$.

If $\fs \neq \mathfrak{i}, \mathfrak{j(\lambda)}$ then $\fs$ is a
generic point in $\fB(\SL(2))$ and we apply Theorem~\ref{one}.
Part~(1) of the conjecture is proved.

%
%

Recall that a representation of $G$ is called \emph{elliptic} if
its character is not identically zero on the elliptic set of $G$.

It follows from
\cite[Theorem~3.4]{Go} (see also \cite[Prop.~1]{SM}) that the (normalized)
induced representation $\Ind_{TU}^G(\lambda \otimes 1)$ has an elliptic
constituent (and then the other constituent is also elliptic) if and only
if the character $\chi$ is of order $2$.

We then have
$$\Ind_{TU}^G(\lambda \otimes 1)=\pi^+\oplus\pi^-,$$
and we have the identity
$\theta^+=\theta^-=0$, between the characters $\theta^+$, $\theta^-$ of
$\pi^+$, $\pi^-$.

Let $s \in\Cset$, $|s| = 1$. The corresponding (normalized) induced
representation will be denoted $\pi(s)$:
\[ \pi(s): = \Ind_{TU}^G(\chi_s \lambda \otimes 1).\]
The representations $\pi(1)$, $\pi(-1)$ are reducible, and split
into irreducible components: \[\pi(1) = \pi^+(1) \oplus \pi^-(1)\]
\[\pi(-1) = \pi^+(-1) \oplus \pi^-(-1)\] These are
\emph{elliptic} representations: their characters \[\theta^+(1),
\theta^-(1), \theta^+(-1), \theta^-(-1)\]are not identically zero
on the elliptic set, although we do have the identities \[
\theta^+(1) + \theta^-(1) = 0\]
\[\theta^+(-1) + \theta^-(-1) = 0.\]


Concerning infinitesimal characters, we have
\[inf.ch.\, \pi^+(1) = inf.ch.\, \pi^-(1) = \lambda\]
\[inf.ch.\, \pi^+(-1) = inf.ch.\, \pi^-(-1) = (-1)^{\val_F} \otimes
\lambda.\]

Note that $E_{\sigma} = \mathbb{T}$, and recall that $\pi(s) \cong
\pi(s^{-1})$. The induced bijection
\[
\Prim^t \cH^{\fs}(G) \to \widetilde{E_{\sigma}}/W_{\ft}\] is as
follows. With $s \in \Cset, |s| = 1$:

\[
\pi(s) \mapsto \{s,s^{-1}\}, \quad \quad s^2 \neq 1\]
\[\pi^+(1) \mapsto 1,\quad \pi^-(1) \mapsto a\]
\[\pi^+(-1) \mapsto -1,\quad \pi^-(-1) \mapsto b\]
where $a,b$ are the two isolated points in the compact extended
quotient of $\mathbb{T}$ by $\Zset/2\Zset$.  We note that the pair
$\pi^+(1), \pi^-(1)$ are \emph{L-indistinguishable}.  The
corresponding Langlands parameter fails to distinguish them from
each other.  The pair $\pi^+(-1), \pi^-(-1)$ are also
L-indistinguishable.

The unramified unitary principal series is defined as follows:
\[ \omega(s): = \Ind_{TU}^G \,(\chi_s \otimes 1).\]  Recall that $\omega(s) = \omega(s^{-1})$.
The representation
$\omega(-1)$ is reducible:
\[
\omega(-1) = \omega^+(-1) \oplus \omega^-(-1).\]

For that part of the tempered spectrum which admits non-zero
Iwahori-fixed vectors, the induced bijection
\[
\Prim^t \cH^{\fs}(G) \to \widetilde{E_{\sigma}}/W_{\ft}\] is as
follows. With $s \in \Cset, |s| = 1$:
\[
\omega(s) \mapsto \{s,s^{-1}\}, \quad \quad s \neq -1\]
\[\omega^+(-1) \mapsto -1, \quad \omega^-(-1) \mapsto c\]
\[\St(2) \mapsto d\]
where $c,d$ are the two isolated points in the compact extended
quotient of $\mathbb{T}$ by $\Zset/2\Zset$, and $\St(2)$ denotes
the Steinberg representation of $\SL(2)$.

These maps are induced by the geometric equivalences and are
therefore not quite canonical, because the geometric equivalences
are not canonical.

%

If $\fs \neq \mathfrak{i}, \mathfrak{j(\lambda)}$ then $\fs$ is a
generic point in $\fB(\SL(2))$ and the induced bijection is as
follows:
\[\Prim^t \,\cH^{\fs}(G) \to \widetilde{E_{\sigma}}/W_{\ft}\]
takes the form of the identity map $\mathbb{T} \to \mathbb{T}$ or
the identity map $pt \to pt$.
\end{proof}

There are $2$ non-generic points in $\fB(\SL(\mathbb{Q}_p))$ with
$p > 3$; there are $4$ non-generic points in
$\fB(\SL(2,\mathbb{Q}_2))$.


\section{Iwahori-Hecke algebras}
\label{general features}

The proof of Theorem~\ref{one} shows that $\cH^\fs(G)$ is Morita
equivalent to a unital $k$-algebra which we will denote by
$\cA_\fs$.
The next step in proving the
conjecture will be to relate this algebra $\cA_\fs$ to a {\it
generalized Iwahori-Hecke algebra}, as defined below.

\par
\smallskip

Let $W'$ be a Coxeter group with generators $(s)_{s\in S}$ and relations
\[(ss')^{m_{s,s'}}=1,\quad\text{for any $s,s'\in S$ such that
$m_{s,s'}<+\infty$,}\] and let $L$ be a weight function on $W'$,
that is, a map $L\colon W'\to\Zset$ such that $L(ww')=L(w)+L(w')$
for any $w$, $w'$ in $W'$ such that $\ell(ww')=\ell(w)+\ell(w')$, where $\ell$
is the usual length function on $W'$. Clearly, the function $\ell$
is itself a weight function.

Let $\Omega$ be a finite group acting on the Coxeter system
$(W',S)$. The group $W:=W'\rtimes\Omega$ will be called an
\emph{extended Coxeter group}. We extend $L$ to $W$ by setting
$L(w\omega):=L(w)$, for $w\in W'$, $\omega\in\Omega$.

Let $A:=\Zset[v,v^{-1}]$ where $v$ is an indeterminate. We set
$u:=v^2$ and $v_s:=v^{L(s)}$ for any $s\in S$. Let $\bar{}\colon A\to A$ be
the ring involution which takes $v^n$ to $v^{-n}$ for any $n\in\Zset$.

Let $\cH(W,u)=\cH(W,L,u)$ denote the $A$-algebra defined
by the generators $(T_s)_{s\in S}$ and the relations
\[(T_s-v_s)(T_s+v_s^{-1})=0\quad\text{for $s\in S$,}\]
\[\underbrace{T_sT_{s'}T_s\cdots}_{\text{$m_{s,s'}$ factors}}=
\underbrace{T_{s'}T_sT_{s'}\cdots}_{\text{$m_{s,s'}$ factors}},
\quad \text{ for any $s\ne s'$ in $S$ such that
$m_{s,s'}<+\infty$.}\]
For $w\in W$, we define $T_w\in\cH(W,u)$ by $T_w=T_{s_1}T_{s_2}\cdots
T_{s_m}$, where $w=s_1s_2\cdots s_m$ is a reduced expression in $W$.
We have $T_1=1$, the unit element of
$\cH(W,u)$, and $(T_w)_{w\in W}$ is an $A$-basis of $\cH(W,u)$. The $v_s$ are
called the parameters of $\cH(W,u)$.

Let $\cH(W',u)$ be the $A$-subspace of $\cH(W,u)$ spanned by all $T_w$
with $w\in W'$.
For each $q\in \Cset^\times$, we set $\cH(W,q):=\cH(W,u)\otimes_{A}\Cset$,
where $\Cset$ is regarded as an $A$-algebra with $u$ acting as
scalar multiplication by $q$.
The algebras of the form $\cH(W,q)$ where $W$ is an extended Coxeter
group and $q\in\Cset$ will be called \emph{extended Iwahori-Hecke algebras}.
In the case when the Coxeter group $W'$ is an affine Weyl group, we will
say that $\cH(W,q)$ is an \emph{extended affine Iwahori-Hecke algebra}.

\par
\medskip

We now observe that $W_\ft$ is a (finite) extended Coxeter
group. Indeed, there exists a root system $\Phi_\ft$ with
associate Weyl group denoted $W'_\ft$ and a subset $\Phi_\ft^+$ of
positive roots in $\Phi_\ft$, such that, setting
\[C_\ft:=\left\{w\in W_\ft\,:\,w(\Phi^+_\ft)\subset\Phi_\ft^+\right\},\]
we have
\[W_\ft=W'_\ft\rtimes C_\ft.\]
This follows from \cite[Prop.~4.2]{Hei} and \cite[Lem.~2]{Ho}.

\par
\medskip

It is expected, and proved, using the theory of types of
\cite{BKtyp}, for level-zero representations in \cite{M},
\cite{M2}, for principal series representations of split groups in
\cite{Rps}, for the group $\GL(n,F)$ in \cite{BK}, \cite{BKss},
for the group $\SL(n,F)$ in \cite{GR2}, for the group $\Sp(4)$ in
\cite{BB2}, and for a large class of representations of classical
groups in \cite{Ki}, \cite{Ki2}, that there exists always an
extended affine Iwahori-Hecke algebra $\cH_\fs'$ such that the
following holds:
\begin{enumerate}
\item there exists a (finite) Iwahori-Hecke algebra $H'_\fs$ with
corresponding Coxeter group $W'_\ft$ and a Laurent polynomial
algebra $\cB_\ft$ satisfying
$\cH_\fs'=H'_\fs\otimes_\Cset\cB_\ft$; \item there exists a
two-cocycle $\mu\colon C_\ft\times C_\ft\to\Cset^\times$ and an
injective homomorphism of groups $\iota\colon C_\ft\to
\Aut_{\Cset-\alg}\cH_\fs'$ such that $\cA_\fs$ is Morita
equivalent to the twisted tensor product algebra \linebreak
$\cH_\fs'\tilde\otimes_\iota\Cset[C_\ft]_\mu$.
\end{enumerate}

In the case of $\GL(n,F)$ (see \cite{BK}), and in the case of
principal series representations of split groups with connected
centre (see \cite{Rps}), we always have $C_\ft=\{1\}$. The
references quoted above give examples in which $C_\ft\ne\{1\}$.
The results in \cite{GR2} also show that the algebra
$\cH_\fs'\tilde\otimes_\iota\Cset[C_\ft]_\mu$ is not always
isomorphic to an extended Iwahori-Hecke algebra.

There are no known example in which the cocycle $\mu$ is non-trivial. In
the case of unipotent level zero representations \cite{Lpadic},
\cite{LpadicII},
of principal series representations \cite{Rps}, and of the group $\Sp_4$
\cite{BB2}, it has been proved that $\mu$ is trivial.

\smallskip
From now we restrict attention to the case where $C_\ft=\{1\}$, so that
$\cA_\fs$ is expected to be Morita equivalent to a generalized
affine Iwahori-Hecke algebra $\cH_\fs'$.

In particular, if $L$ is a torus and $\fs=[T,1]_G$, then $\cA_\fs$ is isomorphic
to the commuting algebra $\cH(G//I)$ in $G$ of the induced representation from
the trivial representation of an Iwahori subgroup $I$ of $G$. We have
\[\cH(G//I)\simeq \cH(W,q_F),\] where $q_F$ is the order
of the residue field of $F$ and $W$ is defined as follows (see for
instance \cite[\S 3.2 and 3.5]{Ca}). Here the weight function
is taken to be equal to the length function. In particular, we are
in the equal parameters case.

Let $T$ be a maximal split torus in $G$, and let $X^*(T)$,
$X_*(T)$ denote its groups of characters and cocharacters,
respectively. Let $\Phi(G,T)\subset X^*(T)$,
$\Phi^\vee(G,T)\subset X_*(T)$ be the corresponding root and
coroot systems, and $W_f$ the associated (finite) Weyl group. Then
\[W=X_*(T)\rtimes W_f,\]
Now let $X'_*(T)$ denote the subgroup of $X_*(T)$ generated by $\Phi^\vee(G,T)$.
Then
$W':=X'_*(T)\rtimes W_f$ is a Coxeter group (an affine Weyl group) and
$W=W'\rtimes\Omega$, where $\Omega$ is the group of elements in $W$ of
length zero.

Let $^LG^0$ be the Langlands dual group of $G$, and let $^LT^0$ denote
the Langlands dual of $T$, a maximal torus of $^LG^0$.
By Langlands duality, we have
\[
W = X_*(T) \rtimes W_f = X^*(^LT^0) \rtimes W_f.\]
 The isomorphism
\[ ^LT^0 \cong \Psi(T), \quad  t \mapsto \chi_t\] is fixed by the relation
\[
\chi_t(\phi(\varpi_F)) = \phi(t)\] for $t \in\,  ^LT^0, \phi \in
X_*(T) = X^*(^LT^0)$, and $\varpi_F$  a uniformizer in $F$. This
isomorphism commutes with the $W_f$-action, see \cite[Section
I.2.3]{GeSh}.

The group $W_f$ acts on $^LT^0$, and we form the quotient variety $^LT^0/W_f$.

Let $\mathfrak{i} \in \fB(G)$ be determined by the cuspidal pair
$(T,1)$. We have
\[ \fZ^{\mathfrak{i}} = \mathbb{C}[^LT^0/W_f],\]
\[\widetilde{\fZ}^{\mathfrak{i}} = \mathbb{C}[\widetilde{^LT^0}/W_f].\]

\section{The asymptotic Hecke algebra} \label{algJ}

There is a unique algebra involution $h\mapsto h^{\dag}$ of $\cH(W',u)$
such that $T_s^{\dag}=-T_s^{-1}$ for any $s\in S$, and
a unique endomorphism $h\mapsto\bar{h}$ of $\cH(W',u)$ which
is $A$-semilinear with respect to $\bar{}\colon A\to A$ and
satisfies $\bar{T}_s=T_s^{-1}$ for any $s\in S$.
Let \[A_{\le
0}:=\bigoplus_{m\le 0}\Zset\, v^m
=\Zset[v^{-1}],\;\;\;
A_{<0}:=\bigoplus_{m<0}\Zset\, v^m,\]
\[\cH(W',u)_{\le 0}:=\bigoplus_{w\in W'} A_{\le 0}\,T_w,\;\;
\cH(W',u)_{<0}:=\bigoplus_{w\in W'} A_{<0}\,T_w.\]
Let $z\in W'$. There is a unique $c_z\in \cH(W',u)_{\le 0}$ such that
$\bar c_z=c_z$ and $c_z=T_z\mod \cH(W',u)_{<0}$, \cite[Theorem~5.2~(a)]{Lbook}.
We write $c_z=\sum_{y\in W'}p_{y,z}\,T_y$, where $p_{y,z}\in A_{\le 0}$.
For $y\in W'$, $\omega,\omega'\in\Omega$, we define $p_{y\omega,z\omega'}$
as $p_{y,z}$ if $\omega=\omega'$ and as $0$ otherwise. For $w\in W$, we set
$c_w:=\sum_{y\in W}p_{y,w}\,T_y$.
Then it follows from \cite[Theorem~5.2~(b)]{Lbook} that
$(c_w)_{w\in W}$ is an $A$-basis of $\cH(W,u)$.

For $x$, $y$, $z$ in $W$, we define $f_{x,y,z}\in A$ by
\[T_xT_y=\sum_{z\in W} f_{x,y,z}\,T_z.\]
From now on, we assume that $W'$ is a bounded weighted Coxeter
group, that is, that there exists
an integer $N\in\Nset$ such that
$v^{-N}f_{x,y,z}\in A_{\le 0}$ for all $x$, $y$, $z$ in $W$.

For $x,y,z$ in $W$, we define $h_{x,y,z}\in A$ by
\[c_x\cdot c_y=\sum_{z\in W}h_{x,y,z}\,c_z.\]
It follows from  \cite[\S 13.6]{Lbook} that, for any $z\in W$, there exists an
integer $\ba(z)\in[0,N]$ such that
\[h_{x,y,z}\in v^{\ba(z)}\,\Zset[v^{-1}]\quad \text{for all $x,y\in W$,}\]
\[h_{x,y,z}\notin v^{\ba(z)-1}\,\Zset[v^{-1}]\quad \text{for some $x,y\in W$.}\]
Let $\gamma_{x,y,z^{-1}}$ be the coefficient of $v^{\ba(z)}$ in
$h_{x,y,z}$.

Let $\Delta(z)\ge 0$ be the integer defined by
\[p_{1,z}=n_zv^{-\Delta(z)}\,+\,\text{strictly smaller powers of
$v$},\;\;\; n_z\in\Zset-\{0\},\]
and let $\cD$ denote the following
(finite) subset of $W$:
\[\cD:=\left\{z\in W\,:\,\ba(z)=\Delta(z)\right\}.\]

In \cite[chap. 14.2]{Lbook} Lusztig stated a list of 15
conjectures $P_1$, $\ldots$, $P_{15}$ and proved them in several
cases \cite[chap. 15, 16, 17]{Lbook}. Assuming the validity of the
conjectures, Lusztig was able to define in \cite{Lbook} partitions
of $W$ into left cells, right cells and two-sided cells, which
extend the theory of Kazhdan-Lusztig from the case of equal
parameters (that is, $v_s=v_{s'}$ for any $(s,s')\in S^2$) to the
general case. In the case of equal parameters the conjectures
mentioned above are known to be true.

From now on we shall assume the validity of these conjectures. Let
us recall some of them below:
\begin{enumerate}
\item[{\bf P1.}] For any $z\in W$ we have $\ba(z)\le\Delta(z)$.
\item[{\bf P2.}] If $d\in\cD$ and $x,y\in W$ satisfy
$\gamma_{x,y,d}\ne 0$, then $x=y^{-1}$. \item[{\bf P3.}] If $y\in
W$, there exists a unique $d\in\cD$ such that
$\gamma_{y,y^{-1},d}\ne 0$. \item[{\bf P4.}] If $z'\le_{\cL\cR}z$
then $\ba(z')\ge\ba(z)$. Hence, if $z'\sim_{\cL\cR}z$, the
$\ba(z')=\ba(z)$. \item[{\bf P5.}] If $d\in\cD$, $y\in W$,
$\gamma_{y^{-1},y,d}\ne 0$, then $\gamma_{y^{-1},y,d}=n_d=\pm 1$.
\item[{\bf P6.}] If $d\in \cD$, then $d^2=1$. \item[{\bf P7.}] For
any $x$, $y$, $z$ in $W$ we have $\gamma_{x,y,z}=\gamma_{y,z,x}$.
\item[{\bf P8.}] Let $x$, $y$, $z$ in $W$ be such that
$\gamma_{x,y,z}\ne 0$. Then $x\sim_{\cL}y^{-1}$,
$y\sim_{\cL}z^{-1}$, $z\sim_{\cL}x^{-1}$.
\end{enumerate}

For any $z\in W$, we set $\hat n_z:=n_d$ where $d$ is the unique
element of $\cD$ such that $d\sim_{\cL}z^{-1}$.

\smallskip

Let $J$ denote the free Abelian group with basis $(t_w)_{w\in W}$.
We set
\[t_x\cdot t_y:=\sum_{z\in W}\gamma_{x,y,z^{-1}}\,t_z.\]
(This is a finite sum.) This defines an associative ring structure
on $J$. The ring $J$ is called the \emph{based ring} of $W$. It
has a unit element $\sum_{d\in\cD}t_d$ (see \cite[\S 18.3]{Lbook}).

The $\Cset$-algebra $J(W):=J\otimes_\Zset\Cset$ is called the
\emph{asymptotic Hecke algebra} of $W$.

\par
\smallskip
According to property (P8), for each two-sided cell $\bc$ in $W$, the subspace
$J_\bc$ spanned by the $t_w$, $w\in\bc$, is a two-sided ideal of $J$.
The ideal $J_\bc$ is in fact an associative ring with unit
$\sum_{d\in\cD\cap\bc}T_d$, which is called the based ring of the
two-sided cell $\bc$, \cite[\S 18.3]{Lbook}, and
\begin{equation} \label{decJ}
J=\bigoplus_\bc J_\bc
\end{equation}
is a direct sum decomposition of $J$ as a ring.

\smallskip

Let $J(W,u):=A\otimes_\Zset J$. We recall from
\cite[Theorem~18.9]{Lbook} (which extends \cite[2.4]{LCellsII})
that the $A$-linear map $\phi\colon \cH(W,u)\to J(W,u)$ given by
\[\phi(c_w^{\dag}):=\sum_{z\in W,\, d\in\cD\atop
\ba(d)=\ba(z)}h_{x,d,z}\,\hat n_z\,t_z\quad (x\in W)\] is a
homomorphism of $A$-algebra with unit (note that the conjecture ($P_{15}$)
is used here).

\smallskip

Let
\begin{equation} \label{phiq}
\phi_q\colon \cH(W,q)\to
J\otimes_\Zset\Cset\end{equation}
be the $\Cset$-algebra homomorphism induced by $\phi$.

Let $\cH(W,q)^{\ge i}$ be the $\Cset$-subspace of $\cH(W,q)$ spanned by
all the $c^{\dag}_w$ with $w\in W$ and $\ba(w)\ge i$. This a
two-sided ideal of $\cH(W,q)$, because of \cite[\S 13.1]{Lbook} and (P7).
Let \[\cH(W,q)^i:=\cH(W,q)^{\ge i}/\cH(W,q)^{\ge
i+1};\] this is an $\cH(W,q)$-bimodule. It has as $\Cset$-basis the
images $[c^{\dag}_w]$ of the $c^{\dag}_w\in \cH(W,q)^{\ge i}$ such
that $\ba(w)=i$.

We may regard $\cH(W,q)^i$ as a $J$-bimodule with multiplication
defined by the rule:
\begin{align}
t_x\,*\,[c^{\dag}_w]&=\sum_{z\in W\atop
\ba(z)=i}\gamma_{x,w,z^{-1}}\,\hat n_w\,\hat n_z\,c^{\dag}_z,\cr
[c^{\dag}_w]\,*\,t_x&=\sum_{z\in W\atop
\ba(z)=i}\gamma_{w,x,z^{-1}}\,\hat n_w\,\hat n_z\,c^{\dag}_z,\quad
(w,x\in W, \ba(w)=i).
\end{align}
We have (see \cite[18.10]{Lbook}):
\begin{equation} \label{hf}
hf\,=\,\phi_q(h)\,*\,f,\;\;\text{for all $f\in \cH(W,q)^i$, $h\in
\cH(W,q)$.}
\end{equation}
On the other side:
\begin{equation} \label{jfh}
(j*f)h=j*(fh),\;\;\text{for all $f\in\cH(W,q)^i$, $h\in\cH(W,q)$, $j\in J$.}
\end{equation}
Let \[f_i:=\sum_{d\in\cD\atop \ba(d)=i}[c_d^{\dag}]\,\in\cH(W,q)^i.\]
\begin{lem} \label{txfi}
We have
\[t_x\,*\,f_i\,=\,f_i\,*\,t_x\,=\,\hat n_x\,[c_x^{\dag}].\]
\end{lem}
\begin{proof}
By definition of $*$, we have
\[t_x\,*\,f_i=\sum_{d\in\cD\atop \ba(d)=i}\sum_{z\in W\atop
\ba(z)=i}\gamma_{x,d,z^{-1}}\,\hat n_d\,\hat n_z\,c^{\dag}_z.\]
Since (because of P7) $\gamma_{x,d,z^{-1}}=\gamma_{d,z^{-1},x}=
\gamma_{z^{-1},x,d}$, it follows from (P2) that if
$\gamma_{x,d,z^{-1}}\ne 0$ then $z=x$. Then, using (P3) and (P5),
we obtain $t_x\,*\,f_i=\hat n_x\,c^{\dag}_x$.
The proof for $f_i\,*\,t_x$ is similar.
\end{proof}

Let $k$ denote the centre of $\cH(W,q)$.
Lusztig proved the following result in \cite[Proposition~1.6~(i)]{LCellsIII}
in the case of equal parameters.
Our proof will follow the same lines.

\begin{prop}
The centre of $J\otimes_\Zset\Cset$ contains $\phi_q(k)$.
\end{prop}
\begin{proof}
It is enough to show that
$\phi_q(z)\,\cdot \,t_x=t_x\,\cdot\,\phi_q(z)$ for any
$z\in k$, $x\in W$. Assume that $\ba(x)=i$. Let $z\in k$.  Using
Lemma~\ref{txfi}, we obtain
\[(\phi_q(z)t_x)*f_i=\phi_q(z)*t_x*f_i=\hat n_x\phi_q(z)*[c^{\dag}_x]=
\hat n_x z[c^{\dag}_x].\] On the other side, using
equation~(\ref{hf}), we get
\[(t_x\phi_q(z))*f_i=t_x*(\phi_q(z)*f_i)=t_x*(zf_i).\]
Since $zf_i=f_iz$, it gives
\[(t_x\phi_q(z))*f_i=t_x*(f_iz),\]
and then, using equation~(\ref{jfh}) and again Lemma~\ref{txfi},
we obtain
\[(t_x\phi_q(z))*f_i=(t_x*f_i)z=\hat
n_x [c_x^{\dag}]z.\] Now, since $z\in k$, we have
$z[c^{\dag}_x]=[c_x^{\dag}]z$. Hence
\begin{equation} \label{interm}
(\phi_q(z)t_x)*f_i=(t_x\phi_q(z))*f_i.
\end{equation}
It follows from the combination of (P4) and (P8) that
$\gamma_{x,y,z}\ne 0$ implies $\ba(x)=\ba(y)=\ba(z)$. Hence we
have
\[\phi_q(z)t_x=\sum_{x'\in W\atop \ba(x')=i}\alpha_{x'}t_{x'},\]
\[t_x\phi_q(z)=\sum_{x'\in W\atop \ba(x')=i}\beta_{x'}t_{x'},\]
with $\alpha_{x'}$, $\beta_{x'}$ in $\Cset$. Then~(\ref{interm})
implies that
\[\sum_{x'\in W\atop \ba(x')=i}\alpha_{x'}[c^{\dag}_{x'}]=
\sum_{x'\in W\atop \ba(x')=i}\beta_{x'}[c^{\dag}_{x'}].\] Hence
$\alpha_{x'}=\beta_{x'}$ for all $x'\in W$ such that $\ba(x')=i$.
It gives $\phi_q(z)t_x=t_x\phi_q(z)$, as required.
\end{proof}

\par
\smallskip

\begin{rem} {\rm
The above proposition provides $J\otimes_\Zset\Cset$ (and also
each $J_\bc$) with a structure of $k$-algebra.  This $k$-algebra
structure is not canonical: it depends on $q$.  Our move (3)
precisely allows us to pass from one $k$-algebra structure,
depending on $q_1$, to another $k$-algebra structure, depending on
$q_2$.}
\end{rem}

\par
\smallskip
From now on we will assume that the weight function is equal
to the length function $\ell$. We will assume that $q=1$, in which
case $H(W,q)$ is the group algebra of $W$; or $q$ is not a root of
unity, in which case we can take for $q$ the order $q_F$ of the
residue field of $F$.

Let $E$ be a simple $\cH(W,q)$-module (resp.
$J\otimes_\Zset\Cset$-module). We attach to $E$ an integer $\ba_E$
by the following two requirements:
\begin{enumerate}
\item
$c_wE=0$ (resp. $t_wE=0$) for any $w$ with $\ba(w)>\ba_E$;
\item
$c_wE\ne 0$ (resp. $t_wE\ne 0$) for some $w$ such $\ba(w)=\ba_E$.
\end{enumerate}

Then Lusztig proved in \cite[Cor.~3.6]{LCellsIII} (see also
\cite[Th.~8.1]{LAst}) that there is a unique bijection $E\mapsto
E'$ between the set of isomorphism classes of simple
$\cH(W,q)$-modules and the set of isomorphism classes of simple
$J\otimes_\Zset\Cset$-modules such that $\ba_{E'}=\ba_E$ and such
that the restriction of $E'$ to $\cH(W,q)$ via $\phi_q$ is an
$\cH(W,q)$-module with exactly one composition factor isomorphic
to $E$ and all other composition factors of the form $\bar E$ with
$\ba_{\bar E}<\ba_E$.

As shown in \cite[Th.~9]{BN}, it follows that $\phi_q$ is spectrum
preserving with respect to filtrations.
Hence
\begin{equation} \label{HtoJ}
\cH(W,q)\asymp J\otimes_\Zset\Cset.
\end{equation}

\par
\medskip

Let $G$ be a connected $F$-split adjoint simple $p$-adic group.
By
Langlands duality we have
\begin{equation} \label{leW}
W: =  X_*(T) \rtimes W_f =
X^*(^LT^0(\mathbb{C})) \rtimes W_f.\end{equation}
Lusztig proved in \cite[Theorem~4.8]{LCellsIV} that the unipotent conjugacy
classes in $^LG^0$ are in bijection with the two-sided cells in $W$.

Let $J$ be the based ring attached to $W$. We fix a two-sided cell
$\bc$ in $W$. Let $\cO_\bc$ be the unipotent conjugacy class in
$^LG^0$ corresponding to $\bc$. Let
$\varphi\colon\SL(2)(\Cset)\to{}^LG^0$ be a homomorphism of
algebraic groups such that $u=\varphi\left(\begin{matrix} 1&1\cr
0&1\end{matrix}\right)$ belongs to $\cO_\bc$ and let $F_\bc$ be
a maximal reductive algebraic subgroup of the centralizer
$\Cent_{{}^LG^0}(u)$.  The reductive group $F_\bc$ may be
disconnected: the identity component of $F_\bc$ will be denoted
$F_{\bc}^0$.

Let $Y$ be a finite $F_\bc$-set (that is, a set with an algebraic action
of $F_\bc$; thus, $F_{\bc}^0$ acts trivially). An $F_\bc$-vector bundle on
$Y$ is a collection of finite dimensional $\Cset$-vector spaces $V_y$
($y\in Y$) with a given algebraic representation of $F_\bc$ on
$\bigoplus_{y\in Y}V_y$ such that $G\cdot V_y=V_{gy}$ for all $g\in
F_\bc$, $y\in Y$. We now consider the finite $F_\bc$-set $Y\times Y$ with
diagonal action of $F_c$ and denote by $K_{F_\bc}(Y\times Y)$ the Grothendieck
group of the category of $F_\bc$-vector bundles on $Y\times Y$. One can
define an associative ring structure on $K_{F_\bc}(Y\times Y)$ (see
\cite[\S 10.2]{LCellsIV}).

Then the conjecture of Lusztig in \cite[\S 10.5]{LCellsIV} states in
particular that there should exist a finite $F_\bc$-set $Y$ and a
bijection $\pi$ from $\bc$ onto the set of irreducible $F_\bc$-vector
bundles on $Y\times Y$ (up to isomorphism) such the $\Cset$-linear map
$J_\bc\to K_{F_\bc}(Y\times Y)\otimes\Cset$ sending $t_w$ to $\pi(w)$ is
an algebra isomorphism (preserving the unit element).

Let $|Y|$ denote the cardinality of $Y$. This number is expected
to be the number of left cells contained in $\bc$. When $F_\bc$ is
connected, $K_{F_\bc}(Y\times Y)$ is isomorphic to the
$|Y|\times|Y|$ matrix algebra
$\Mat_{|Y|}(R_{\mathbb{C}}({F_\bc}))$ over the (complexified)
rational representation ring $R_{\mathbb{C}}({F_\bc})$ of $F_\bc$.
It is important to note: when $F_\bc$ is connected, the Lusztig
conjecture asserts that $J_\bc$ is Morita equivalent to a
\emph{commutative} algebra.

The Lusztig conjecture has been proved by Xi for any two-sided
cell $\bc$ when $G$ is one of the following groups $\GL(n)$,
$\PGL(n)$, $\SL(2)$, $\SO(5)$ and $G_2$, and for the lowest
two-sided cell $\bc_0$ (see next section) when $G$ is any
connected $F$-split adjoint simple $p$-adic group.

\section{The ideal $J_{\bc_0}$ in $J$}

As above we assume that $G$ is a connected $F$-split adjoint simple
$p$-adic group.
Let $J$ be the based ring attached to $W$, with $W$ as in~(\ref{leW}).
The centralizer of $1$ is of course $^LG^0$.
Under the bijection cited above, the unipotent
class $1$ corresponds to the \emph{lowest two-sided cell} $\bc_0$, that is
the subset of all the elements $w$ in $W$ such that $\ba(w)$ equals the number
of positive roots in the root system of $W_f$.

Xi proved the Lusztig conjecture for
this ideal $J_{\bc_0}$ in \cite[Theorem~1.10]{Xic0}.
According to his result, we have a ring isomorphism
\[J_{\bc_0} \cong \Mat_{|W_f|}(R_{\mathbb{C}}(^LG^0)).\]

The character map $Ch$ creates an isomorphism
\[R_{\mathbb{C}}(^LG^0) \cong (R_{\mathbb{C}}(^LT^0))^{W_f}.\]
The $W_f$-invariant subring of the (complexified) representation
ring of $^LT^0$ is precisely the coordinate ring of the quotient
torus $^LT^0/W_f$.

Since
\[\Psi(T) =  {^L}T^0\]
we have a Morita equivalence
\[
J_{\bc_0} \sim \fZ^{\mathfrak{i}}(G)\] where $\mathfrak{i}$ is the
quotient variety $\Psi(T)/W_f$. Therefore, we obtain the following
result.

\begin{thm} \label{Jczero}
Let $G$ be a connected $F$-split adjoint simple $p$-adic group. There is
a Morita equivalence between $J_{\bc_0}$ and the coordinate ring
of the Bernstein variety $\Psi(T)/W_f$.
\end{thm}

According to our conjecture, \emph{the other ideals $J_{\bc}$
account (up to geometric equivalence) for the rest of the extended
quotient of $\Psi(T)$ by $W_f$}.

The classical Satake isomorphism is an isomorphism between the
spherical Hecke algebra $\cH(G//K)$ and the ring
$R_{\Cset}(^LG^0)$.  Further, a theorem of Bernstein (see e.g.
\cite[Proposition 8.6]{L0}) asserts that the centre $Z(\cH(G//I))$
of the Iwahori-Hecke algebra $\cH(G//I)$ is also isomorphic to
$R_{\Cset}(^LG^0)$.

At this point, we need the map $\phi_{q_F,\bc_0}$ defined in section
1.7 of Xi's paper \cite{Xic0}. This map is the composition of
$\phi_{q_F}$ and of the projection of $J$ onto $J_{\bc_0}$.

Xi has proved in \cite[Theorem 3.6]{Xic0} that the image
$\phi_{q_F,\bc_0}(Z(H(G//I))$ is the centre $Z(J_{\bc_0})$ of the
algebra $J_{\bc_0}$. This creates the following diagram:

\[
\begin{CD}
\cH(G//K) @>>> R_{\Cset}(^LG^0)\\
@VVV                @VVV\\
\cH(G//I) @>>\phi_{q_F,\bc_0}> J_{\bc_0} \end{CD}\] in which the top
horizontal map is the Satake isomorphism, the left vertical map is
induced by the inclusion $K \subset I$, the right vertical map
sends $R_{\Cset}(^LG^0)$ onto the centre of $J_{\bc_0}$ and the
bottom horizontal map is Xi's map $\phi_{q_F,\bc_0}$.  The vertical
maps are injective.  We expect that this diagram is commutative.

\section{The Hecke algebra of $\GL(n)$}

\begin{thm} \label{GL}  The conjecture is true for $\GL(n)$.
\end{thm}

\begin{proof} In this proof, we follow \cite{BP} rather closely;
we have refined the proof at certain points.  The occurrence of an
extended quotient in the smooth dual of $\GL(n)$ was first
recorded in \cite{HP}, in the context of Deligne-Langlands
parameters.

Let $G:=\GL(n)$, $\fs=[L,\sigma]_G\in\fB(G)$ and
$\ft=[L,\sigma]_L\in\fB(L)$. We can think of $\ft$ as a vector of
irreducible supercuspidal representations of smaller general
linear groups. If the vector is
\[(\sigma_1, \ldots, \sigma_1, \ldots, \sigma_t, \ldots,
\sigma_t)\] with $\sigma_i$ repeated $\evarepsilon_i$ times, $1 \leq i \leq
t$, and $\sigma_1$, $\ldots$, $\sigma_t$ pairwise distinct (after
unramified twist) then we say that $\ft$ has {\it exponents}
$\evarepsilon_1$, $\ldots$, $\evarepsilon_t$.

Each representation $\sigma_i$ of $G_i:=\GL(m_i)$ has a {\it
torsion number}: the order of the cyclic group of all those
unramified characters $\eta$ for which $\sigma_i \otimes \eta
\cong \sigma_i$.  The torsion number of $\sigma_i$ will be denoted
$r_i$.

Hence
\[L\simeq\prod_{i=1}^t G_i^{\evarepsilon_i}\;\;\text{ and }\;\;
\sigma\simeq\bigotimes_{i=1}^t\sigma_i^{\otimes \evarepsilon_i},\] Each
$\sigma_i$ contains a maximal simple type $(K_i,\lambda_i)$ in
$G_i$ \cite{BK}. Let
\[K_L:=\prod_{i=1}^t K_i^{\evarepsilon_i} \;\;\text{ and
}\;\;\tau_L:=\bigotimes_{i=1}^t\lambda_i^{\otimes \evarepsilon_i}.\] Then
$(K_L,\tau_L)$ is a $\ft$-type in $L$.
We have
\[W_\ft\simeq\prod_{i=1}^tS_{\evarepsilon_i}.\]
Let $W_{\evarepsilon_i}$ denote the extended affine Weyl group associated to
$\GL(\evarepsilon_i,\Cset)$.

Let $(K,\tau)$ be a semisimple $\fs$-type, see \cite{BK, BKtyp,
BKss}.  It is worth pointing out that we do not need the type
explicitly. Instead, we need certain items attached to the type:
the idempotent $e_{\tau}$ and the endomorphism-valued Hecke
algebra $\cH(G,\tau)$.   Let $e_{\tau}$ be the idempotent attached
to the type $(K,\tau)$ as in \cite[Definition 2.9]{BKtyp}:

\[
e_{\tau}(x) = \begin{cases}(\text{vol}\, K)^{-1}(\dim\,\tau)\,tr
(\tau(x^{-1}))&\text{ if $x \in K$,}\cr 0&\text{ if $x \in G, x
\notin K$.}\end{cases}\]

 The idempotent $e_{\tau}$ is
then a \emph{special} idempotent in the Hecke algebra $\cH(G)$
according to \cite[Definition 3.11]{BKtyp}. Let $\cH = \cH(G)$.
It follows from \cite[\S 3]{BKtyp} that
\[
\cH^{\fs}(G) = \cH * e_{\tau} * \cH.\]

 We then have a Morita
equivalence
\[
\cH * e_{\tau} * \cH \sim_{morita} e_{\tau} * \cH * e_{\tau}.\]
 Now let $\mathcal{H}(K, \tau)$ be the
endomorphism-valued Hecke algebra attached to the semisimple type
$(K, \tau)$.   By \cite[2.12]{BKtyp} we have a canonical
isomorphism of unital $\mathbb{C}$-algebras :

\[\mathcal{H}(G, \tau) \otimes_{\mathbb{C}} \End_{\mathbb{C}}W
\cong e_{\tau} * \mathcal{H}(G) * e_{\tau}\] so that the algebra
$e_{\tau}
* \mathcal{H}(G) * e_{\tau}$ is Morita equivalent to the
algebra $\mathcal{H}(G, \tau)$.   Now we quote the main theorem
for semisimple types in $\GL(n)$ \cite[1.5]{BKss}:  there is an
isomorphism of unital $\mathbb{C}$-algebras

\[\mathcal{H}(G, \tau) \cong   \bigotimes_{i=1}^t
\cH(W_{\evarepsilon_i},q_F^{r_i}).\]

The factors $\mathcal{H}(W_{\evarepsilon_i}, q_F^{r_i})$ are (extended) affine
Hecke algebras whose structure is given explicitly in \cite[5.4.6,
5.6.6]{BK}.

We conclude that
\[
\cH^{\fs}(G) \asymp \bigotimes_{i=1}^t \cH(W_{\evarepsilon_i},q_F^{r_i}).\]

On the other hand, from~(\ref{HtoJ}), we have
\[\bigotimes_{i=1}^t \cH(W_{\evarepsilon_i},q_F^{r_i})
 \asymp\bigotimes_{i=1}^t J(W_{\evarepsilon_i}).\]

Finally we will prove that that
\[\bigotimes_{i=1}^t J(W_{\evarepsilon_i}) \asymp \widetilde{\fZ^{\fs}}.\]

Let $^LT^0$ be the maximal standard torus of $^LG^0=\GL(n, \Cset)$ and
let $W$ be the extended affine Weyl group associated to $\GL(n, \Cset)$.
We have
$W:=X^*(^LT^0)\rtimes S_{n}=W_n$. For each two-sided cell
$\bc$ of $W$ we have a corresponding partition $\lambda$ of $n$.
Let $\mu$ be the dual partition of $\lambda$. Let $u$ be a unipotent element
in $\GL(n,\Cset)$ whose Jordan blocks are determined
by the partition $\mu$. Let the distinct parts of the dual
partition $\mu$ be $\mu_1, \dots, \mu_p$ with $\mu_r$ repeated
$n_r$ times, $1 \leq r \leq p$.

Let $C_G(u)$ be the centralizer of $u$ in $G = \GL(n , \Cset)$.
Then the maximal reductive subgroup $F_{\mathbf{c}}$ of $C_G(u)$
is isomorphic to $\GL(n_1, \bbc) \times \GL(n_2, \bbc)\times \dots
\times \GL(n_p, \Cset)$. For the non-trivial combinatorics which
underlies this statement, see \cite[\S 2.6]{MG}.

Let $J$ be the based ring of $W$.
For each two-sided cell in $W$, let $|Y|$ be the number of left cells
contained in $\mathbf{c}$. The Lusztig conjecture says that there
is a ring isomorphism
$$
J_{\mathbf{c}} \simeq \Mat_{|Y|} (R_{F_{\mathbf{c}}}), \quad t_w
\mapsto \pi(w)
$$
where $R_{F_{\mathbf{c}}}$ is the rational representation ring of
$F_{\mathbf{c}}$. This conjecture for $\GL(n, \Cset)$ has been
proved by Xi \cite[1.5, 4.1, 8.2]{XM}.

Since $F_{\mathbf{c}}$ is isomorphic to a direct product of the
general linear groups $\GL(n_i, \Cset)$ ($1 \leq i \leq p$) we see
that $R_{F_{\mathbf{c}}}$ is isomorphic to the tensor product over
$\Zset$ of the representation rings $R_{\GL(n_i, \Cset)}$, $1\leq
i \leq p$. For the ring $R({\GL(n, \Cset))}$ we have
\[R(\GL(n,\Cset)) = \mathbb{Z}[X^*(T(\Cset))]^{S_n}\] where $T(\Cset)$ is the standard
maximal torus in $\GL(n,\Cset)$, and $X^*(T(\Cset))$ is the set of
rational characters of $T(\Cset)$, by \cite[Chapter VIII]{Bou}.
Therefore we have

$$
R_{F_{\mathbf{c}}} \otimes_\Zset \Cset \simeq \Cset
[\text{Sym}^{n_1}\Cset^\times \times \cdots\times \text{Sym}^{n_p}
\Cset^\times].$$

Let $\gamma \in S_n$ have cycle type $\mu$, let $X =
(\bbc^\times)^n$. Then
$$
\begin{array}{rcl}
X^\gamma & \simeq & (\Cset^\times)^{n_1} \times \dots \times (\Cset^\times )^{n_p}\\
Z(\gamma) & \simeq & (\Zset/\mu_1 \Zset) \wr S_{n_1} \times \dots
\times
(\Zset/\mu_p \Zset ) \wr S_{n_p}\\
X^\gamma /Z(\gamma) & \simeq & \text{Sym}^{n_1} \Cset^\times
\times \dots \times \text{Sym}^{n_p}\Cset^\times
\end{array}
$$
and so
$$
R_{F_{\mathbf{c}}}\otimes _\Zset \Cset \simeq \Cset [X^\gamma
/Z(\gamma)].
$$
Then, using~(\ref{decJ}), we obtain
$$
J\otimes_\Zset \Cset = \bigoplus_{\mathbf{c}} (J_{\mathbf{c}} \otimes
_\Zset \Cset)  \sim \bigoplus_{\mathbf{c}} (R_{F_{\mathbf{c}}}
\otimes_\Zset \Cset )
 \simeq \Cset [\widetilde{X}/S_n].
$$
The algebra $J \otimes_\Zset \Cset$ is Morita equivalent to a
reduced, finitely generated, commutative unital $\Cset$-algebra,
namely the coordinate ring of the extended quotient
$\widetilde{X}/S_n$.   This finishes the proof of part (1) of the
conjecture.

Part (2) of the conjecture for $\GL(n)$ is a consequence of
\cite[Theorem Theorem 5.1]{P}.
\end{proof}

\section{The Iwahori ideal in $\cH(\PGL(n))$}

Let $G = \PGL(n)$, let $T$ be its standard maximal torus. Let $W:
= X_*(T) \rtimes W_f$.   Then $^LG^0 = \SL(n,\mathbb{C})$ is the
Langlands dual group. Its maximal torus will be denoted $^LT^0$.

The discrete group $W$ is an \emph{extended Coxeter group}: \[ W =
\langle s_1, s_2, \ldots, s_n \rangle \rtimes
\mathbb{Z}/n\mathbb{Z}\] where $\Zset/n\Zset$ permutes cyclically
the generators $s_1, \ldots, s_n$. We have
\[
\cH(W,q_F) = \cH(G//I).\]

The symmetric group $W_f = S_n$ acts on $^LT^0$ by permuting
coordinates, and we form the quotient variety $^LT^0/S_n$.

Let $\mathfrak{i} \in \fB(G)$ be determined by the cuspidal pair
$(T,1)$. We have
\[ \fZ^{\mathfrak{i}} = \mathbb{C}[^LT^0/S_n],\]
\[\widetilde{\fZ}^{\mathfrak{i}} = \mathbb{C}[\widetilde{^LT^0}/S_n].\]

\begin{thm} \label{PGL}
Let $\cH^{\mathfrak{i}}(G)$ denote the Iwahori ideal in $\cH(G)$.
Then $\cH^{\mathfrak{i}}(G)$ is geometrically equivalent to the
extended quotient of $^LT^0$ by the symmetric group $S_n$:
\[
\cH^{\mathfrak{i}}(G) \asymp \mathbb{C}[\widetilde{^LT^0}/S_n].\]
\end{thm}

\begin{proof}  The non-unital algebra $\cH^{\mathfrak{i}}(G)$ is Morita
equivalent to the unital affine Hecke algebra $\cH(W,q_F)$:
\[
\cH^{\mathfrak{i}}(G) = \cH e \cH \sim_{morita} e \cH e \cong
\cH(W,q_F).\]

From~(\ref{HtoJ}), we have
\[\cH(W,q_F)\asymp J\otimes_\Zset\Cset.\]
The Langlands dual of $\PGL(n,F)$ is $\SL(n,\Cset)$.
For each two-sided cell
$\mathbf{c}$ of $W$ we have a corresponding partition $\lambda$ of $n$.
Let $\mu$ be the dual partition of $\lambda$. Let $u$ be a unipotent element
in $\SL(n,\Cset)$ whose Jordan blocks are determined
by the partition $\mu$. Let the distinct parts of the dual
partition $\mu$ be $\mu_1 < \cdots < \mu_p$ with $\mu_r$ repeated
$n_r$ times, $1 \leq r \leq p$.

Let $C_G(u)$ be the centralizer of $u$ in $G = \SL(n , \Cset)$.
Then the maximal reductive subgroup $F'_{\lambda}$ of $C_G(u)$ is
isomorphic to $(\GL(n_1, \bbc) \times \GL(n_2, \bbc)\times \dots
\times \GL(n_p, \Cset)) \cap \SL(n,\mathbb{C})$.  For details of
the injective map \[(\GL(n_1, \bbc) \times \GL(n_2, \bbc)\times
\dots \times \GL(n_p, \Cset)) \cap \SL(n,\mathbb{C})
\longrightarrow \SL(n,\mathbb{C})\] see \cite{MG}.

As a special case, let the two-sided cell $\bc$ correspond to the
partition $\lambda = (1,1,1, \ldots,1)$ of $n$. Then the dual
partition $\mu = (n)$. The unipotent matrix $u$ has one Jordan
block, and its centralizer $C_G(u) = Z$ the centre of $\SL(n,
\Cset)$. The maximal reductive subgroup $F'_{\lambda}$ of $C_G(u)$
is the finite group $Z$. This is the case $p= 1$, $\mu_1 = n$, $n_1 =
1$.

By the theorem of Xi \cite[8.4]{XM} we have
\[
J_{\bc}\otimes_{\mathbb{Z}} \mathbb{C} \sim_{morita}
R_{\mathbb{C}}(F'_{\lambda}) = R_{\mathbb{C}}(Z) = \mathbb{C}^n.\]

Let $\gamma$ have cycle type $(n)$.  Then the fixed set
$(^LT^0)^{\gamma}$ comprises the $n$ fixed points
\[
\diag(\omega^j, \ldots, \omega^j) \in ^LT^0\] where
$\omega = \exp(2\pi i/n)$ and $0 \leq j \leq n-1$. These $n$ fixed
points correspond to the $n$ generators in the commutative ring
$\mathbb{C}^n$.

We expect that the corresponding points in $\Irr(\PGL(N))$ arise
as follows. The unramified unitary twist
\[ z^{\val \circ \det} \otimes \St(n)\] of the Steinberg
representation of $\GL(n)$ has trivial central character if and
only if $z$ is an $n$th root of unity.   For these values $1,
\omega, \omega^2, \ldots, \omega^{n-1}$ of $z$, we obtain $n$
irreducible smooth  representations of $\PGL(n)$.

From now on, we will assume that $\lambda \neq (1,1,1,\ldots,1)$.
Then $F'_{\lambda}$ is a \emph{connected} Lie group.

We will write $T_{\lambda}(\Cset)$ for the standard maximal torus
of $F'_{\lambda}$. The Weyl group is then
\[W(\lambda) = S_{n_1} \times \cdots \times S_{n_p}.\]

According to Bourbaki \cite[Chapter 8]{Bou}, the map $Ch$, sending
each (virtual) representation to its (virtual) character, creates
an isomorphism:
\[Ch: R(F'_{\lambda}) \cong \mathbb{Z}[X^*(T_{\lambda}(\Cset)]^{W(\lambda)}.\]

Note that a complex linear combination of rational characters of
$T_{\lambda}(\Cset)$ is precisely a regular function on
$T_{\lambda}(\Cset)$.

For each two-sided cell $\mathbf{c}$ of $W$ the $\Zset$-submodule
$J_{\mathbf{c}}$ of $J$, spanned by all $t_w$, $w \in \mathbf{c}$,
is a two-sided ideal of $J$. The ring $J_{\mathbf{c}}$ is the
based ring of the two-sided cell $\mathbf{c}$. Now apply the
theorem of Xi \cite[8.4]{XM}. We get
\[
J_{\bc}\otimes_{\Zset} \Cset \sim_{morita} R_{\Cset}(F'_{\lambda})
\cong \Cset[T_{\lambda}(\Cset)]^{W(\lambda)} \cong
\Cset[T_{\lambda}(\Cset)/W(\lambda)].\]

Let $\gamma \in S_n$ have cycle type $\mu$. Then the
$\gamma$-centralizer is a direct product of wreath products:
\[
Z(\gamma)  \simeq  (\Zset/\mu_1 \Zset) \wr S_{n_1} \times \dots
\times (\Zset/\mu_p \Zset ) \wr S_{n_p}.\]

The image of $T_{\lambda}(\Cset)$ in the inclusion
$T_{\lambda}(\Cset) \to ^LT^0$ is precisely the subtorus of
$^LT^0$ fixed by  $\Zset/\mu_1 \Zset  \times \dots \times
\Zset/\mu_p \Zset$.  We therefore have
\[ (^LT^0)^\gamma /Z(\gamma) \simeq \, T_{\lambda}(\Cset)/W(\lambda).\]
We conclude that
$$
R(F'_{\lambda})\otimes _\Zset \Cset \simeq \Cset [(^LT^0)^\gamma
/Z(\gamma)]
$$
Then
$$
J\otimes_\Zset \Cset = \oplus_{\mathbf{c}} (J_{\mathbf{c}} \otimes
_\Zset \Cset)  \sim \oplus_{\mathbf{c}} (R(F'_{\lambda})
\otimes_\Zset \Cset )
 \simeq \Cset [\widetilde{^LT^0}/S_n]
$$
The algebra $J \otimes_\Zset \Cset$ is Morita equivalent to a
reduced, finitely generated, commutative unital $\Cset$-algebra,
namely the coordinate ring of the extended quotient
$\widetilde{^LT^0}/S_n$.
\end{proof}

\section{The Iwahori ideal in $\cH(\SO(5))$}

Let $G$ denote the special orthogonal group $\SO(5,F)$. We view it
as the group of elements of determinant $1$ which stabilise the
symmetric bilinear form
$$\left(\begin{matrix}
0&0&0&0&1\\
0&0&0&1&0\\
0&0&1&0&0\\
0&1&0&0&0\\
1&0&0&0&0\end{matrix}\right).$$ Let $T$ be the group of diagonal
matrices
$$\left(\begin{matrix}\lambda_1&0&0&0&0\\
0&\lambda_2&0&0&0\\
0&0&1&0&0\\
0&0&0&\lambda_2^{-1}&0\\
0&0&0&0&\lambda_1^{-1}\end{matrix}\right),\quad
\text{$\lambda_1,\lambda_2\in F^\times$.}$$ The extended affine
Weyl group $W=X_*(T)\ltimes W_f$ is of type $\tilde B_2$, with
$W_f\simeq S_2\ltimes (\Zset/2)^2$ a finite Weyl group of type
$B_2$.

Here $^LG^0=\Sp(4,\Cset)$, and $^L T^0$ is the group of diagonal
matrices
$$\diagt(t_1,t_2):=\left(\begin{matrix}t_1&0&0&0\\
0&t_2&0&0\\
0&0&t_2^{-1}&0\\
0&0&0&t_1^{-1}\end{matrix}\right),\quad \text{$t_1,t_2\in
\Cset^\times$.}$$ We have $\Nor_{^LG^0}(^L T^0)/^L T^0\simeq W_f$.
The group $^LG^0 = \Sp(4,\Cset)$ is simply connected. In \cite[\S
11.1]{XL}, Xi has proved Lusztig's conjecture for the group $W$.

The extended Coxeter group $W=W'\rtimes\Omega$ has four two-sided
cells $c_\id$, $c_1$, $c_2$ and $c_0$ (see \cite[\S 11.1]{XL}):
\[c_\id=\left\{w\in W \,:\, \ba(w)=0\right\}=\{\id,\omega\}=\Omega,\]
\[c_1=\left\{w\in W \,:\, \ba(w)=1\right\},\]
\[c_2=\left\{w\in W \,:\, \ba(w)=2\right\},\]
\[c_0=\left\{w\in W \,:\, \ba(w)=4\right\} \;\;\text{(the lowest two-sided
cell)}.\] We have \[ J = J_{\bc_\id} \oplus J_{\bc_1} \oplus
J_{\bc_2} \oplus J_{\bc_0}.\]

Let $\mathfrak{i} \in \fB(G)$ be determined by the cuspidal pair
$(T,1)$. We have
\[ \fZ^{\mathfrak{i}} = \mathbb{C}[^LT^0/W_f],\quad
\widetilde{\fZ}^{\mathfrak{i}} =
\mathbb{C}[\widetilde{^LT^0}/W_f].\]

\begin{lem} \label{thirdlemma}
Let \[\mathbb{L}: = \mathbb{C}[t,t^{-1}]^{\Zset/2\Zset}\] denote the balanced
Laurent polynomials in one indeterminate $t$, where the generator
$\alpha$ of $\Zset/2\Zset$ acts as follows: $\alpha(t) = t^{-1}$.
Then the coordinate algebra of the extended quotient of $^LT^0$ by
$W_f$
 is the $\Cset$-algebra
\[
\Cset^5 \oplus \mathbb{L}^3 \oplus \Cset[^LT^0/W_f].\]
\end{lem}
\begin{proof}  Let $X:={}^LT^0$.
The $8$ elements $\gamma_1$, $\ldots$, $\gamma_8$ of $W_f$ can be
described as follows:
\[\begin{matrix}
\gamma_1(\diagt(t_1,t_2))=(\diagt(t_1,t_2)) &
\gamma_2(\diagt(t_1,t_2))=\diagt(t_2,t_1)\\
\gamma_3(\diagt(t_1,t_2))=\diagt(t_1^{-1},t_2)&
\gamma_4(\diagt(t_1,t_2))=\diagt(t_1,t_2^{-1})\\
\gamma_5(\diagt(t_1,t_2))=\diagt(t_2,t_1^{-1})&
\gamma_6(\diagt(t_1,t_2))=\diagt(t_1^{-1},t_2^{-1})\\
\gamma_7(\diagt(t_1,t_2))=\diagt(t_2^{-1},t_1)&
\gamma_8(\diagt(t_1,t_2))=\diagt(t_2^{-1},t_1^{-1})\\
\end{matrix}.\]
We have $\gamma_5=\gamma_2\gamma_3=\gamma_4\gamma_2$,
$\gamma_6=\gamma_5^2=\gamma_3\gamma_4=\gamma_4\gamma_3$,
$\gamma_7=\gamma_5^3=\gamma_2\gamma_4=\gamma_3\gamma_2$,
$\gamma_8=\gamma_4\gamma_2\gamma_4=\gamma_2\gamma_4\gamma_3$. The
elements $\gamma_2$, $\gamma_3$, $\gamma_4$, $\gamma_6$ and
$\gamma_8$ are of order $2$, the elements $\gamma_5$ and
$\gamma_7$ are of order $4$. We obtain
$$X^{\gamma_1}=X, \quad X^{\gamma_2}=\left\{\diagt(t,t)\,:\,t\in\Cset^\times
\right\},$$
$$X^{\gamma_3}=X^{\gamma_4}=\left\{\diagt(1,t_2),\diagt(-1,t_2)\,:\,t_2\in\Cset^\times
\right\},$$
$$X^{\gamma_5}=X^{\gamma_7}=\left\{\diagt(1,1),\diagt(-1,-1)\right\},$$
$$X^{\gamma_6}=\{\diagt(1,1),\diagt(1,-1),\diagt(-1,1),\diagt(-1,-1)\}.$$
The elements $\gamma_1$, $\gamma_6$ are central, and we have
$$Z(\gamma_2)=\left\{\gamma_1,\gamma_2,\gamma_6,\gamma_8\right\},$$
$$
Z(\gamma_3)=\left\{\gamma_1,\gamma_3,\gamma_4,\gamma_6\right\},\quad
Z(\gamma_7)=\left\{\gamma_1,\gamma_5,\gamma_6,\gamma_7\right\}.$$
There are five $W_f$-conjugacy classes:
$$\{\gamma_1\},\{\gamma_6 \},\{\gamma_2,\gamma_8\},\{\gamma_3,\gamma_4\},
\{\gamma_5,\gamma_7\}.$$

As representatives, we will take $\gamma_1, \gamma_2, \gamma_3,
\gamma_5, \gamma_6$.

\begin{itemize}

\item We have \[\Cset[X^{\gamma_6}/Z(\gamma_6)] =
\Cset[X^{\gamma_6}/W_f] = \Cset \oplus \Cset \oplus \Cset\] since
there are three $W_f$-orbits in $X^{\gamma_6}$, namely
\[\{\diagt(1,1)\}, \{\diagt(1,-1), \diagt(-1,1)\}, \{\diagt(-1,-1)\}\]

\item We have \[\Cset[X^{\gamma_5}/Z(\gamma_5)] =
\Cset[X^{\gamma_5}] = \Cset \oplus \Cset = R_{\Cset}(Z) =
\Cset[\Omega] = J_{\bc_{\id}}\] since $Z(\gamma_5)$ acts trivially
on $X^{\gamma_5}$.

\item We have \[X^{\gamma_3} = \{\diagt(1,t) ; t \in
\Cset^{\times}\} \sqcup \{\diagt(-1,t): t \in \Cset^{\times}\} \]
The $Z(\gamma_3)$-orbit of $\diagt(1,t)$ is the unordered pair
$\{\diagt(1,t),\diagt(1,t^{-1})\}$ and the $Z(\gamma_3)$-orbit of
$\diagt(-1,t)$ is the unordered pair
$\{\diagt(-1,t),\diagt(-1,t^{-1})\}$.  Therefore we have
\[\Cset[X^{\gamma_3}/Z(\gamma_3)] = \mathbb{L} \oplus
\mathbb{L}.\]

\item We have $X^{\gamma_2} = \{\diagt(t,t) : t \in
\Cset^{\times}\} \cong \{t: t \in \Cset^{\times}\}$. The
$Z(\gamma_2)$-orbit of the point $t$ is the unordered pair
$\{t,t^{-1}\}$. So we have
\[X^{\gamma_2}/Z(\gamma_2) \cong \Cset^{\times}/\Zset/2\] and
$\Cset[X^{\gamma_2}/Z(\gamma_2)] = \mathbb{L}$.

\item We have $\Cset[X^{\gamma_1}/Z(\gamma_1)] =
\Cset[^LT^0/W_f]$.
\end{itemize}
\end{proof}

The reductive group $F_{\bc_{\id}}$ is the center of $^LG^0$,
$F_{\bc_1}=(\Zset/2\Zset)\ltimes\Cset^\times$ where $\Zset/2\Zset$
acts on $\Cset^\times$ by $z\mapsto z^{-1}$, and
$F_{\bc_2}=(\Zset/2\Zset)\times\SL(2,\Cset)$ and there is a ring
isomorphism $J_{\bc_2}\simeq \Mat_4(R_{F_{\bc_2}})$, where
$R_{F_{\bc_2}}$ is the rational representation ring of
$F_{\bc_2}$, \cite[Theorem~11.2]{XL}.

 We have $F_{\bc_1} = \langle \alpha \rangle \ltimes \Cset^{\times}$
 where $\alpha$ generates
$\Zset/2\Zset$.  Note the crucial relation \[ \alpha z = z^{-1}
\alpha\]with $z \in \Cset^{\times}$. The (semisimple) conjugacy
classes in $F_{\bc_1}$ are:
\[
\{1\}, \quad \{-1\}, \quad \{\{z,z^{-1}\}: z \in \Cset^{\times},
z^2 \neq 1\}, \quad \alpha \cdot \Cset^{\times}\]

\begin{lem} \label{fourthlemma}
Let $\mathbb{M}: = \Cset[t,t^{-1}]$ denote the Laurent polynomials
in one indeterminate $t$.  We have
\[ J_{\bc_1} \asymp \mathbb{C} \oplus (\mathbb{M}\rtimes \Zset/2\Zset)\]
where the generator $\alpha$ of $\Zset/2\Zset$ acts as follows:
$\alpha(t) = t^{-1}$.
\end{lem}
\begin{proof} Let $F = F_{\bc_1}$, $F = \langle \alpha
\rangle \ltimes \mathbb{C}^{\times}$. We have to construct the
simple $J_{\bc_1}$-modules explicitly, following Xi \cite[p. 51,
107]{XL}. We will use Xi's explicit proof of the Lusztig
conjecture for $B_2$. Let $Y = \{1,2,3,4\}$ be the $F$-set such
that as $F$-sets we have $\{1\} \cong \{2\} \cong F/F$ and
$\{3,4\} \cong F/F^0$. The simple $J_{\bc_1}$-modules are given by
\[
E_{s,\rho}: = \Hom_{A(s)}(\rho,H_*(Y^s))\] where $A(s)$ denotes
the component group $C_F(s)/C_F(s)^0$ and $\rho$ is a simple
$A(s)$-module which appears in the homology group $H_*(Y^s)$.  The
set $Y^s$ denotes the $s$-fixed set: $Y^s = \{y \in Y: sy = y\}$.
The pair $(s,\rho)$ is chosen up to $F$-conjugacy.

(A). $s = 1$, $\rho = 1$. Then $Y^s = \{1,2,3,4\}$. Also $ H_*(Y^s)$
is the free $\Cset$-vector space $V: = \Cset^4$ on $\{1,2,3,4\}$:
we will denote its basis by $\{e_1,e_2,e_3,e_4\}$. $A(s) =
\Zset/2\Zset$.  The generator of $A(s)$ permutes $e_3, e_4$. Let
$V_1$ denote the span of $e_1,e_2, e_3 + e_4$, let $V_2$ denote
the span of $e_3 - e_4$.  Then
\[
E_{1,1}: =  \Hom_{A(s)}(\rho, V) = V_1 = \Cset^3.\]

(B). $s = 1$, $\rho = \epsilon$ where $\epsilon$ is the sign
representation of $\Zset/2\Zset$.  We have
\[
E_{1,\epsilon}: =  \Hom_{A(s)}(\rho, V) = V_2 = \Cset.\] Note that
we have
\[
E_{1,1} \oplus E_{1,\epsilon} = V = \Cset^4\] as $A(s)$-modules.

(C). $s = -1$, $\rho = 1$. We have $Y^s = Y$, $A(s) = \Zset/2\Zset$
and
\[
E_{-1,1}: = \Hom_{A(s)}(\rho, V) =  V_1 = \Cset^3\]

(D). $s = -1$, $\rho = \epsilon$. We have $Y^s = Y$, $A(s) =
\Zset/2\Zset$ and
\[
E_{-1,\epsilon}: = \Hom_{A(s)}(\rho, V) =  V_2 = \Cset\] Note that
we have
\[
E_{-1,1} \oplus E_{-1,\epsilon} = V = \Cset^4\] as $A(s)$-modules.

(E) $s = z$, $\rho = 1$ where $z \in \Cset^{\times}$, $z^2 \neq 1$. We
have $Y^s = Y$, $A(s) = \{1\}$ and
\[
E_{z,1}: = \Hom(\Cset, V) =  V = \Cset^4\]

(F). $s = \alpha$, $\rho = 1$. We have $Y^s = \{1,2\}$, $H_*(Y^s) =
\Cset^2$, $A(s) = \{1,-1,\alpha, -\alpha\} = \Zset/2\Zset \times
\Zset/2\Zset$. We have
\[
E_{\alpha,1}: = \Hom_{A(s)}(\Cset, \Cset^2) = \Cset^2.\]

This concludes the list of simple $J_{\bc_1}$-modules. We now turn
to the Lusztig-Xi isomorphism of unital algebras:
\[J_{\bc_1} \cong K_F(Y \times Y) \otimes_{\mathbb{Z}} \mathbb{C}.\]
The equivariant $K$-theory $K_F(Y \times Y)\otimes_{\Zset}\Cset$
is equipped with the convolution product. The action of $F$ on the
set $Y$ leads to the following description. We identify $K_F(Y
\times Y) \otimes_{\Zset} \Cset$ with the $\Cset$-algebra of $4
\times 4$ matrices $(a_{ij})$ where $a_{11}, a_{12}, a_{21},
a_{22} \in R(F)\otimes_{\Zset} \Cset$ and all other entries are in
$R(\mathbb{C}^{\times})\otimes_{\Zset} \Cset$, subject to the
following conditions:
\[a_{14} = \overline{a_{13}},\, a_{24} = \overline{a_{23}},\, a_{41} =
\overline{a_{31}},\, a_{42} = \overline{a_{32}},\, a_{44} =
\overline{a_{33}},\, a_{43}= \overline{a_{34}}\] where, for all $z
\in \Cset^{\times}$,
\[ \overline{a_{ij}}(z) = a_{ij}(z^{-1}).\]

Let
\[\mathbb{M}: = \mathbb{C}[t,t^{-1}]\]
denote the Laurent polynomials in one indeterminate $t$. We have
an injective homomorphism of unital $\mathbb{C}$-algebras:
\[\psi : \mathbb{C} \oplus (\mathbb{M} \rtimes \mathbb{Z}/2\mathbb{Z})
\longrightarrow K_F(Y \times Y) \otimes_{\Zset}\Cset\]

\[(\lambda,p + \alpha \cdot q) \mapsto
\left(\begin{matrix}\lambda &0&0&0\\
0& \lambda &0&0\\
0&0&p&\overline{q}\\
0&0&q&\overline{p}\end{matrix}\right)\]

We claim that this map is a spectrum-preserving morphism of unital
finite-type $k$-algebras.   We view each entry in $K_F(Y \times
Y)$ as a virtual character of $F$ or $\Cset^{\times}$.  We can
then \emph{evaluate} each matrix entry at $z \in \Cset^{\times}$.
Evaluation at $z \in \mathbb{C}^{\times}$  determines an algebra
homomorphism
\[K_F(Y \times Y) \otimes_{\mathbb{Z}} \mathbb{C} \longrightarrow M_4(\Cset).\]
This homomorphism gives $\Cset^4$ the structure of
$J_{\bc_1}$-module.

\begin{itemize}
\item When $z^2 \neq 1$, this homomorphism gives $\Cset^4$ the
structure of a simple $J_{\bc_1}$-module, namely $E_{z,1}$. \item
When $z = 1$, the module $\Cset^4$ splits into simple
$J_{\bc_1}$-modules of dimensions $3$ and $1$, namely $E_{1,1}$
and $E_{1,\epsilon}$. \item When $z = -1$, the module $\Cset^4$
splits into simple $J_{\bc_1}$-modules of dimensions $3$ and $1$,
namely $E_{-1,1}$ and $E_{-1,\epsilon}$.
\end{itemize}

When restricted to the lower right $2 \times2$ block of the $4
\times 4$ matrix algebra $(a_{ij})$,  each simple $4$-dimensional
module splits into the direct sum of a $2$-dimensional $0$-module
and a $2$-dimensional simple module over the crossed product
$\mathbb{M} \rtimes \Zset/2\Zset$. At $1$ and $-1$, the module
$\Cset^4$ restricts to a $2$-dimensional $0$-module and two
$1$-dimensional simple modules. It now follows that
$\Prim(J_{\bc_1})$ with one deleted point is in bijection (via our
morphism of algebras) with $\Prim(\mathbb{M} \rtimes
\Zset/2\Zset)$. There is one remaining point in
$\Prim(J_{\bc_1})$, namely $E_{\alpha,1}$. This point in
$\Prim(J_{\bc_1})$ maps, via our morphism of algebras, to the one
remaining primitive ideal $0 \oplus (\mathbb{M} \rtimes
\Zset/2\Zset)$.

We conclude from this that
\[
J_{\bc_1} \asymp \Cset \oplus (\mathbb{M} \rtimes \Zset/2\Zset).\]
\end{proof}

\begin{thm} \label{SO}  Let $\cH^{\mathfrak{i}}(\SO(5))$ denote the Iwahori ideal in $\SO(5)$.  We have
\[
\cH^{\mathfrak{i}}(\SO(5)) \asymp
\widetilde{\fZ}^{\mathfrak{i}}(\SO(5)).
\]
\end{thm}
\begin{proof}
 We note that
\begin{itemize}
\item $J_{\bc_\id} = \Cset[\Omega] = \Cset^2$ \item  $J_{\bc_2}
\asymp R_{\Cset}(F_{\bc_2}) = \mathbb{L}^2$ since
\[R_{\Cset}(F_{\bc_2}) = R_{\Cset}(\Zset/2\Zset \times \SL(2,\Cset)) =
R_{\Cset}(\Zset/2\Zset) \otimes R_{\Cset}(\SL(2,\Cset))  =
\mathbb{L} \oplus \mathbb{L}.\]
\end{itemize}

 By Lemmas~\ref{secondlemma} and~\ref{fourthlemma}, we have
\[J =  J_{\bc_\id} \oplus
J_{\bc_1} \oplus J_{\bc_2} \oplus J_{\bc_0} \asymp \Cset^2 \oplus
(\Cset^3 \oplus \mathbb{L}) \oplus \mathbb{L}^2 \oplus
J_{\bc_0}.\] By Lemma~\ref{thirdlemma}, we have
\[\Cset[\widetilde{^LT^0}/W_f]
= \Cset^5 \oplus \mathbb{L}^3 \oplus \Cset[^LT^0/W_f].\] We
conclude, by Theorem~\ref{Jczero}, that
\[\cH^{\mathfrak{i}}(\SO(5)) \asymp J \asymp
\Cset[\widetilde{X}/W_f] =
\widetilde{\fZ}^{\mathfrak{i}}(\SO(5))\] as required.  This
confirms part (1) of our conjecture for the Iwahori ideal of
$\cH(\SO(5))$.
\end{proof}

It is worth noting that, in this example, there are $4$ two sided
cells and $5$ $W_f$-conjugacy classes.

\section{Consequences of the conjecture}
\textsc{Parametrization of the smooth dual}. In this section we
will suppose that the conjecture is true for the Iwahori ideal
$\cH^{\mathfrak{i}}(G)$.   We then have
\[\cH^{\mathfrak{i}}(G) \asymp
\mathcal{O}(\widetilde{^LT^0}/W_f).\]

We now take the primitive ideal space of each side.  We obtain a
bijection
\[
\Prim\,\cH^{\mathfrak{i}}(G) \longleftrightarrow
\widetilde{^LT^0}/W_f.\]

Now $\Prim\,\cH^{\mathfrak{i}}(G)$ may be identified with the
subset $\Irr_I(G)$ of the smooth dual $\Irr(G)$ which admits
nonzero Iwahori-fixed vectors.   This leads to a parametrization
of $\Irr_I(G)$:
\[
\Irr_I(G) \longleftrightarrow \widetilde{^LT^0}/W_f.\]

This is a \emph{parametrization of $\Irr_I(G)$ by the
$\Cset$-points in a complex affine algebraic variety (with several
components)}. This parametrization is not quite canonical: it
depends on the specific finite sequence of elementary steps (and
filtrations) connecting $\cH^{\mathfrak{i}}(G)$ to
$\mathcal{O}(\widetilde{^LT^0}/W_f)$.

This parametrization will assign to each $\omega \in \Irr_I(G)$ a
pair
\[
(s,\gamma)\] with $s \in ^LT^0$, $\gamma$ a $W_f$-conjugacy class.

More generally, the conjecture leads to a parametrization of
$\Irr(G)$ by the $\Cset$-points in a complex affine locally
algebraic variety (with countably many components).  The
dimensions of the components are less than or equal to the rank of
$G$.

\textsc{Langlands parameters}.   We wish to make a comparison with
Langlands parameters.   We first recall some background material.

Let $\cW_F$ denote the Weil group of $F$, and let
$\cI_F\subset\cW_F$ be the inertia subgroup so that
$\cW_F/\cI_F\simeq \Zset$. Denote the Frobenius generator by
$\Frob$ in order that $\cW_F=\cI_F\rtimes\langle\Frob\rangle$.

By \emph{Langlands parameter} we mean a continuous homomorphism
$$\varphi\colon W_F\times\SL(2,\Cset)\to {}^L G^0 $$
which is rational on $\SL(2,\Cset)$ and such that $\varphi(\Frob)$
is semisimple.

Call $\varphi$ {\it unramified} if it has trivial restriction to
$\cI_F$. The unramified Langlands parameters are parameterized by
pairs $(s,u)$ with $s\in {}^LG^0$ semisimple, $u\in {}^LG^0$
unipotent with $sus^{-1}=u^{q_F}$.

To each Langlands parameter $\varphi$ should correspond a finite
packet $\Pi_\varphi$ of irreducible smooth representations of $G$.
According to the Langlands philosophy, the unramified Langlands
parameters should be those for which the representations in
$\Pi_\varphi$ are unipotent in the sense of \cite{Lpadic}. This
has been proved by Lusztig \cite{Lpadic} when $G$ is the group of
$F$-points of a split adjoint simple algebraic group.

More precisely, let $Z$ be the centre of ${}^LG^0$, let $C(g)$
denote the centralizer in ${}^LG^0$ of $g\in {}^LG^0$, and
$C(s,u):=C(s)\cap C(u)$. We denote by $C(s,u)^0$ the connected
component of $C(s,u)$ and set $A(s,u):=C(s,u)/Z\cdot C(s,u)^0$. It
is proved in \cite{Lpadic} that the isomorphism classes of
unipotent representations of $G$ are naturally in one-to-one
correspondence with the set of triples $(s,u,\rho)$ (modulo the
natural action of ${}^LG^0$) where $s$, $u$ as above and $\rho$ an
irreducible representation (up to isomorphism) of $A(s,u)$.

In the case when we restrict ourselves to representations in
$\Irr_I(G)$ the above classification specializes to the
Kazhdan-Lusztig classification, which provides a proof of a
refined form of a conjecture of Deligne and Langlands: if $G$ is
the group of $F$-points of any split reductive group over $F$, the
set $\Irr_I(G)$ is naturally in bijection with the set of triples
$(s,u,\rho)$ as above (modulo the natural action of ${}^LG^0$)
such that $\rho$ appears in the homology of the variety $\cB_u^s$
of the Borel subgroups of ${}^LG^0$ containing both $s$ and $u$
(see \cite{KL}, \cite{R}).

For $\GL(n)$, the simultaneous centralizer $C(s,u)$ is connected,
and so each Deligne-Langlands parameter is a pair $(s,u)$. There
is a bijection between unipotent classes in $\GL(n,\Cset)$ and
partitions of $n$ (via Jordan canonical form). There is a second
bijection between partitions of $n$ and conjugacy classes in $W_f
= S_n$ the symmetric group.   In this way we obtain a perfect
match
\[
(s,u) \longleftrightarrow (s,\gamma).\]

The Deligne-Langlands parameters can be arranged to form an
extended quotient of the complex torus $\Cset^{\times}$ by the
symmetric group $S_n$. The details of this correspondence were
recorded in \cite{HP}.

\smallskip

\textsc{Irreducibility of induced representations}.  If $W_{\ft}$
acts freely then $\widetilde{D_{\sigma}}/W_{\ft} \cong
D_{\sigma}/W_{\ft}$ and the conjecture in this case predicts
irreducibility of the induced representations $\iota_P^G(\sigma
\otimes \chi)$. This situation is discussed in \cite{R1}, with any
maximal Levi subgroup of $G$, and also with the following (non
maximal) Levi subgroup
\[
L = (\GL(2) \times \GL(2) \times \GL(4)) \cap \SL(8) \] of
$\SL(8)$ and $\sigma$ defined as in \cite[p.127]{R1}.

Anne-Marie Aubert, Institut de Ma\-th\'e\-ma\-ti\-ques de Jussieu,
U.M.R. 7586 du C.N.R.S., Paris, France\\
Email: aubert@math.jussieu.fr

Paul Baum, Pennsylvania State University, Mathematics Department,
University Park, PA 16802, USA\\
Email: baum@math.psu.edu

Roger Plymen, School of Mathematics, University of Manchester,
Manchester M13 9PL, England\\
Email: plymen@manchester.ac.uk

\begin{thebibliography}{99}
\bibitem{BC} P. Baum, A. Connes, The Chern character for discrete
groups, A F\^{e}te of topology, Academic Press, New York, 1988, 163--232.
\bibitem{BN} P. Baum, V. Nistor, Periodic cyclic homology of Iwahori-Hecke
algebras, K-Theory 27 (2002), 329--357.
\bibitem{BNotes} J. Bernstein, Representations of $p$-adic groups, Notes by
K.E. Rumelhart, Harvard University 1992.
\bibitem{B} J. Bernstein, P.~Deligne, D.~Kazhdan, M.-F.~Vigneras,
Re\-pr\'e\-sen\-ta\-tions des grou\-pes r\'e\-duc\-tifs sur un corps local,
Travaux en cours, Hermann, 1984.
\bibitem{Bbook} J. Bernstein and S. Gelbart, An introduction to the Langlands program.
Birkhauser, Basel, 2003.
\bibitem{BB2} L.~Blasco, C.~Blondel, Alg\`ebres de Hecke et s\'eries
principales g\'en\'eralis\'ees de $\SP_4(F)$, Proc. London Math.
Soc. (3) 85 (2002), 659--685.
\bibitem{Bou} N. Bourbaki, Lie groups and Lie algebras, Chapters VI, VIII and
IX, Masson, 1981.
\bibitem{BP} J. Brodzki, R.J. Plymen, Complex structure on the smooth dual of
$\GL(n)$, Documenta Math. 7 (2002), 91--112.
63 (2001), no. 2, 364--386.
\bibitem{BK} C.J. Bushnell, P.C. Kutzko, The admissible dual of
$\GL(n)$ via compact open subgroups, Ann. Math. Study 129,
Princeton Univ. Press 1993.
\bibitem{BKtyp} C.J. Bushnell, P.C. Kutzko, Smooth representations of
reductive $p$-adic groups: Structure theory via types, Proc.
London Math. Soc. 77 (1998), 582--634.
\bibitem{BKss} C.J. Bushnell, P.C. Kutzko, Semisimple types in $\GL_n$,
Comp. Math. 119 (1999), 53--97.
\bibitem{Ca} P. Cartier, Representations of $p$-adic groups: a survey,
Proc. Symp. Pure Math. 33 (1979), part 1, 111--155.
\bibitem{CDN} N. Chifan, S. Dascalescu and C. Nastasescu, Wide Morita contexts,
relative injectivity and equivalence results,  J. Algebra 284
(2005), 705--736.
\bibitem{C} J. Cuntz, Morita invariance in cyclic
homology for nonunital algebras, K-Theory 15 (1998), 301 - 305.
\bibitem{CST} J. Cuntz, G. Skandalis, B. Tsygan, Cyclic homology in noncommutative geometry,
EMS 121, 2004.
\bibitem{EH} D. Eisenbud, J. Harris, The geometry of schemes, Springer Graduate
Text 197, 2001.
\bibitem{MG} M. Geck, An introduction to algebraic geometry and
algebraic groups, Oxford Graduate Texts in Mathematics 10, Oxford
University Press, 2003.
\bibitem{GeSh} S. Gelbart, F. Shahidi, Analytic properties of
automorphic $L$-functions, Perspectives in Math. 6 (1988).
\bibitem{Go} D.~Goldberg, $R$-groups and elliptic representations for
$\SL_n$, Pacific J. of Math. 165 (1994), 77--92.
\bibitem{GR2} D.~Goldberg, A.~Roche, Hecke algebras and $\SL_n$-types,
Proc. London Math. Soc. (3)  90  (2005),  no. 1, 87--131.
\bibitem{Hei} V. Heiermann, D\'ecomposition spectrale et
re\-pr\'e\-sen\-ta\-tions sp\'e\-ci\-a\-les d'un groupe
r\'e\-duc\-tif $p$--adique, Journal de l'Institut Math. Jussieu, 3
(2004), 327--395.
\bibitem{HP} J.E. Hodgins and R.J. Plymen, The representation theory of $p$-adic
$\GL(n)$  and Deligne-Langlands parameters, TRIM 10 (1996), 54--72.
\bibitem{Ho} R.~B.~Howlett, Normalizers of parabolic subgroups of
reflection groups, J. London Math. Soc. (2), 21 (1980), 62--80.
\bibitem{KL} D.~Kazhdan, G.~Lusztig, Proof of the Deligne-Langlands
conjecture for Hecke algebras, Invent. Math. 87 (1987) 153--215.
\bibitem{KNS} D.~Kazhdan, V.~Nistor, P.~Schneider, Hochschild and cyclic
homology of finite type algebras, Sel. Math., New ser. 4 (1998), 321--359.
\bibitem{Ki} J.-L.~Kim, Hecke algebras of classical groups over $p$-adic
fields and supercuspidal representations, Amer. J. Math. 121 (1999),
967--1029.
\bibitem{Ki2} J.-L.~Kim, Hecke algebras of classical groups over $p$-adic
fields II, Compositio Math. 127 (2001), 117--167.
\bibitem{Lem} F. Lemmermeyer, Reciprocity laws from Euler to Eisenstein, Springer Monograph
in Mathematics, Berlin, 2000.
\bibitem{JL} J.-L. Loday, Cyclic homology, Springer-Verlag, Berlin, 1992.
\bibitem{L0} G. Lusztig, Singularities, character formulas and a
$q$-analog of weight multiplicities, Ast\'{e}risque 101 --102
(1983), 208--229.
\bibitem{Lpadic} G. Lusztig, Classification of unipotent representations
of simple p-adic groups, Internat. Math. Res. Notices 11 (1995), 517--589.
\bibitem{LpadicII} G. Lusztig, Classification of unipotent representations
of simple p-adic groups, II, Represent. Theory 6
(2002), 243--289.
\bibitem{LCellsII} G. Lusztig, Cells in affine Weyl groups, II, J. Algebra 109
(1987), 536--548.
\bibitem{LCellsIII} G. Lusztig, Cells in affine Weyl groups, III, J. Fac. Sci.
Univ. Tokyo, Sect. IA, Math, 34 (1987), 223--243.
\bibitem{LCellsIV} G. Lusztig, Cells in affine Weyl groups, IV, J. Fac.
Sci. Univ. Tokyo, Sect. IA, Math, 36 (1989), 297--328.
\bibitem{LAst} G. Lusztig, Representations of affine Hecke algebras,
Ast\'erisque 171-172 (1989), 73-84.
\bibitem{Lpadic} G. Lusztig, Classification of unipotent representations
of simple p-adic groups, Internat. Math. Res. Notices 11 (1995),
517--589.
\bibitem{Lbook} G. Lusztig, Hecke algebras with unequal parameters, CRM
Monograph Series, 18, Amer. Math. Soc. 2003.
\bibitem{Mis} P. Mischenko, Invariant tempered distributions on the reductive
$p$-adic group $\GL_n(F_p)$, C. R. Math. Rep. Acad. Sci. Canada 4
(1982),123--127.
\bibitem{M} L. Morris, Tamely ramified intertwining algebras, Invent. Math.
114 (1993), 1--54.
\bibitem{M2} L. Morris, Level zero $G$-types, Compositio Math. 118 (1999),
135--157.
\bibitem{N} J. Neukirch, Algebraic number theory, Springer-Verlag, 1999.
\bibitem{P}  R.J. Plymen, Reduced $C^*$-algebra of the $p$-adic
group $\GL(n)$ II, J. Functional Analysis 196 (2002) 119--134.
\bibitem{R} M.~Reeder, Isogenies of Hecke algebras and a Langlands
correspondence for ramified principal series representations,
Represent. Theory 6 (2002), 101--126.
\bibitem{R1} A. Roche, Parabolic induction and the Bernstein
decomposition, Compositio Math. 134 (2002), 113--133.
\bibitem{Rps} A. Roche, Types and Hecke algebras for principal
series representations of split reductive $p$-adic groups, Ann. scient.
\'Ec. Norm. Sup. {\bf 31} (1998), 361--413.
\bibitem{JR} J. Rogawski, The nonabelian reciprocity law for local
fields, Notices of the American Math. Society, January (2000), 35--41.
\bibitem{SM} G.~Sanje Mpacko, Types for elliptic non-discrete series
representations of $\SL_N(F)$, $N$ prime and $F$ a $p$-adic field,
J. Lie Theory 5 (1995), 101--127.
\bibitem{SZ} P. Schneider, E.-W. Zink, $K$-types for the tempered
components of a $p$-adic general linear group, J. reine angew.
Math. 517 (1999) 161 --208.
\bibitem{Wald} J.-L. Waldspurger,  La formule de Plancherel d'apr\`es
Harish-Chandra. Journal de l'Institut Math. Jussieu 2 (2003)
235--333.
\bibitem{XM} N. Xi, The based ring of two-sided cells of affine Weyl groups
of type $\widetilde{A}_{n-1}$, Amer. Math. Soc. Memoir 749, 2002.
\bibitem{XL} N. Xi, Representations of affine Hecke algebras, Lecture Notes in
Math. 1587, Springer, 1994.
\bibitem{Xic0} N. Xi, The based ring of the lowest two-sided cell of an
affine Weyl group, J. of Algebra 134 (1990) 356--368.
\end{thebibliography}
\end{document}